\documentclass[twoside,10pt]{article}
\usepackage{graphicx}
\usepackage{subfigure}
\usepackage{url}
\usepackage{color}

\usepackage{amsfonts}
\newcommand\Q{\mathbb Q}
\newcommand\R{\mathbb R}
\newcommand\T{\mathbb T}
\newcommand\Z{\mathbb Z}
\newcommand\ee{\mathrm e}
\newcommand\ii{\mathrm i}

\newcommand\eps\varepsilon
\newcommand\A{\mathcal A}
\newcommand\Cc{\mathcal C}
\newcommand\F{\mathcal F}
\newcommand\I{\mathcal I}
\newcommand\J{\mathcal J}
\newcommand\Lc{\mathcal L}
\newcommand\M{\mathcal M}
\newcommand\Ord{\mathcal O}
\newcommand\Pc{\mathcal P}
\newcommand\Tc{\mathcal T}
\newcommand\W{\mathcal W}
\newcommand\Y{\mathcal Y}
\newcommand\Zc{\mathcal Z}

\newcommand\ds{\displaystyle}
\newcommand\dfrac[2]{\ds\frac{#1}{#2}}

\newcommand\p[1]{\left(#1\right)}
\newcommand\pq[1]{\left[#1\right]}
\newcommand\pp[1]{\left\{#1\right\}}
\newcommand\scprod[2]{\left\langle#1,#2\right\rangle}
\newcommand\abs[1]{\left|#1\right|}
\newcommand\absd[1]{\left|#1\right|_2}
\newcommand\norm[1]{\left\|#1\right\|}

\newcommand\tl{\tilde}
\newcommand\wtl{\widetilde}
\newcommand\wh[1]{\widehat{#1}}
\newcommand\whh[1]{\wh{#1}\,'}
\newcommand\ol[1]{\overline{#1}}
\newcommand\tp{^\top}

\newcommand\df{{\rm d}}
\newcommand\Df{{\rm D}}
\newcommand\Id{{\rm Id}}
\newcommand\rint{\mathop{\rm rint}}
\newcommand\ut{{\rm u}} 
\newcommand\st{{\rm s}} 
\newcommand\loc{{\rm loc}}
\newcommand\diag{\mathop{\rm diag}}
\newcommand\SL{\mathop{\rm SL}}
\newcommand\GL{\mathop{\rm GL}}
\newcommand\beps{\bar\eps}    
\newcommand\ceps{\check\eps}  
\newcommand\bzeta{\bar\zeta}
\newcommand\bg{\ol g}
\newcommand\bbf{\ol f}
\newcommand\bb{\ol b}
\newcommand\bF{\ol F}
\newcommand\bR{\ol R}
\newcommand\bx{\bar x}
\newcommand\bbeta{\bar\beta}
\newcommand\arsinh{\mathop{\rm arsinh}}
\newcommand\Lg{\mathop{\rm Lg}} 

\newtheorem{theorem}{Theorem}
\newtheorem{proposition}[theorem]{Proposition}
\newtheorem{lemma}[theorem]{Lemma}
\newcommand\bremark{\noindent{\bf Remark.}\ \ }
\newcommand\eremark{\medskip}
\newcommand\bremarks{\noindent{\bf Remarks.}\ \bnm}
\newcommand\eremarks{\enm\medskip}
\newcommand\proof{\noindent\textbf{\emph{Proof.}}\quad}
\newcommand\sketchproof{\noindent\textbf{\emph{Sketch of the proof.}}\quad}
\newcommand\proofof[1]{\noindent\textbf{\emph{Proof of #1.}}\quad}
\newcommand\qed{\ \ \null\nolinebreak\hfill$\frame{\large\phantom a}$}
\newcommand\paragr[1]{\noindent\textbf{\emph{#1.}}\quad}

\newcommand\beq{\begin{equation}}
\newcommand\eeq{\end{equation}}
\newcommand\bea{\begin{eqnarray}}
\newcommand\eea{\end{eqnarray}}
\newcommand\bean{\begin{eqnarray*}}
\newcommand\eean{\end{eqnarray*}}
\newcommand\btm{\vspace{-\baselineskip}\begin{itemize}}
\newcommand\etm{\end{itemize}\vspace{-\baselineskip}}
\newcommand\bnm{\vspace{-\baselineskip}\begin{enumerate}}
\newcommand\enm{\end{enumerate}\vspace{-\baselineskip}}

\begin{document}

\title{Exponentially small splitting of separatrices\\
  associated to 3D whiskered tori with cubic frequencies
  \ \footnote{This work has been partially supported by
      the Spanish MINECO/FEDER grants
      MTM2015-65715, PGC2018-098676-B-I00 (the~three authors)
      and MTM2016-80117-P (the author MG),
      the Catalan grants 2017SGR1049 (AD and PG) and 2017SGR1374 (MG),
      the Russian Scientific Foundation grant 14-41-00044 (AD and MG),
      and the Juan de la Cierva--Formaci\'on/Incorporaci\'on fellowships
      FJCI-2014-21229 and IJCI-2016-29071 (MG).}}

\author{\sc
  Amadeu Delshams$\,^{1\,,\,2}$,
  \ Marina Gonchenko$\,^1$,
\\[4pt]\sc
  \ Pere Guti\'errez$\,^1$
\\[12pt]
  {\small
  $^1\;$\parbox[t]{5cm}{
    Dep. de Matem\`atiques\\
    Univ. Polit\`ecnica de Catalunya\\
    Av. Diagonal 647, 08028 Barcelona
  }
  \ $^2\;$\parbox[t]{6cm}{
    Lab of Geometry and Dynamical Systems\\
    Univ. Polit\`ecnica de Catalunya\\
    Av. Dr. Mara\~n\'on 44--50, 08028 Barcelona
  }}
\\[25pt]
  {\footnotesize E-mail: \parbox[t]{3.7cm}{
    \texttt{amadeu.delshams@upc.edu}\\
    \texttt{marina.gonchenko@upc.edu}\\
    \texttt{pere.gutierrez@upc.edu}
  }}}
\maketitle

\begin{abstract}
We study the splitting of invariant manifolds of whiskered (hyperbolic) tori
with three frequencies in a nearly-integrable Hamiltonian system, whose
hyperbolic part is given by a pendulum. We consider a 3-dimensional torus with
a fast frequency vector $\omega/\sqrt\varepsilon$, with
$\omega=(1,\Omega,\widetilde\Omega)$ where $\Omega$ is a cubic irrational
number whose two conjugates are complex, and the components of $\omega$
generate the field $\mathbb Q(\Omega)$. A paradigmatic case is the cubic golden
vector, given by the (real) number $\Omega$ satisfying $\Omega^3=1-\Omega$, and
$\widetilde\Omega=\Omega^2$. For such 3-dimensional frequency vectors, the
standard theory of continued fractions cannot be applied, so we develop a
methodology for determining the behavior of the small divisors
$\langle k,\omega\rangle$, $k\in{\mathbb Z}^3$. Applying the
Poincar\'e-Melnikov method, this allows us to carry out a careful study of the
dominant harmonic (which depends on $\varepsilon$) of the Melnikov function,
obtaining an asymptotic estimate for the maximal splitting distance, which is
exponentially small in $\varepsilon$, and valid for all sufficiently small
values of~$\varepsilon$. This estimate behaves like
$\exp\{-h_1(\varepsilon)/\varepsilon^{1/6}\}$ and we provide, for the first
time in a system with 3 frequencies, an accurate description of the (positive)
function $h_1(\varepsilon)$ in the numerator of the exponent, showing that it
can be explicitly constructed from the resonance properties of the frequency
vector $\omega$, and proving that it is a quasiperiodic function (and not
periodic) with respect to $\ln\varepsilon$. In this way, we emphasize the
strong dependence of the estimates for the splitting on the arithmetic
properties of the frequencies.
\par\vspace{12pt}
\noindent\emph{Keywords}:
  splitting of separatrices,
  transverse homoclinic orbits,
  Melnikov integrals,
  cubic frequency vectors
\par\vspace{12pt}
\noindent\emph{AMS subject classification}:
  37J40, 70H08
\end{abstract}

\section{Introduction and setup}

\subsection{Background and state of the art}

In nearly-integrable Hamiltonian systems with $n\ge2$ degrees of freedom,
irregular motion may take place near ($n-1$)-dimensional whiskered tori
(invariant hyperbolic tori) and their whiskers (invariant manifolds).
In adequate scaled canonical coordinates
(see for instance \cite{DelshamsG01,Lochak90,DelshamsGG14a}
and references therein for more details about this introductory paragraph),
these whiskered tori have frequency vectors with fast frequencies and
their non-small hyperbolic part is typically given by a pendulum.
The fundamental phenomenon guaranteeing irregular behavior near these
whiskered tori is the non-coincidence of their whiskers,
which is called the \emph{splitting of separatrices}.
The size of this splitting provides a measure of the irregular motion
(and also of the global instability for $n\ge3$ degrees of freedom)
but is non-easily computable, since it turns out to be exponentially small
with respect to the perturbation parameter.
To worse things, for $n\ge3$, the exponent in the splitting depends strongly
on the arithmetic properties of the ($n-1$)-dimensional frequency vectors
of the whiskered torus.
Fortunately, for $n=3$ the standard theory of continued fractions can be
successfully applied to the $2$-dimensional frequency vectors of the whiskered
tori to compute the splitting.
Nevertheless, for $n\ge4$ degrees of freedom,
the standard theory of continued fractions cannot be applied to
($n-1$)-dimensional frequency vectors, and so far there are no computations of
the exponentially small splitting of separatrices for whiskered tori with
dimension greater or equal than three.

This paper is dedicated to the study and computation of the exponentially small
splitting of separatrices,
in a perturbed Hamiltonian system with 4~degrees of freedom,
associated to a 3-dimensional whiskered torus
with a \emph{cubic frequency vector}.
More precisely, we start with an integrable Hamiltonian $H_0$ possessing
whiskered tori with a \emph{homoclinic whisker} or \emph{separatrix},
formed by coincident stable and unstable whiskers,
and we focus our attention on a concrete torus
with a frequency vector of \emph{fast frequencies}:
\beq\label{eq:omega_eps}
  \omega_\eps = \frac\omega{\sqrt{\eps}}\;,
  \qquad
  \omega=(1,\Omega,\wtl\Omega),
\eeq
with a small (positive) parameter $\eps$, and we assume that
the frequency ratios $\Omega=\omega_2/\omega_1$ and
$\wtl\Omega=\omega_3/\omega_1$ (it can be assumed that $\omega_1=1$)
generate a \emph{complex cubic field}
(also called a \emph{non-totally real} cubic field).
This amounts to assume that $\Omega$ is a \emph{cubic irrational number}
(a real root of a polynomial of degree~3 with rational coefficients,
that is not rational or quadratic) whose two \emph{conjugates are not real},
and $\wtl\Omega=a_0+a_1\Omega+a_2\Omega^2$, with $a_0,a_1,a_2\in\Q$, $a_2\ne0$
(see Section~\ref{sect:resonantseq1} for more details).
A paradigmatic example is the vector $\omega=(1,\Omega,\Omega^2)$,
where $\Omega$ is the \emph{cubic golden number}
(the real number satisfying $\Omega^3=1-\Omega$,
see Section~\ref{sect:cubicgolden}).

If we consider a perturbed Hamiltonian $H=H_0+\mu H_1$, where $\mu$ is small,
in general the whiskers do not coincide anymore.
This phenomenon has got the name of \emph{splitting of separatrices},
which is related to the non-integrability of the system
and the existence of chaotic dynamics,
and plays a key role in the description of Arnold diffusion.
If we assume, for the two involved parameters,
a relation of the form $\mu=\eps^r$
for some $r>0$, we have a problem of singular perturbation
and in this case the splitting is \emph{exponentially small}
with respect to $\eps$.
Our aim is to provide an \emph{asymptotic estimate} for the
\emph{maximal splitting distance}, and to show the dependence of such estimate
on the \emph{arithmetic properties} of the cubic number~$\Omega$.

To provide a measure for the splitting, we can restrict ourselves to
a transverse section to the unperturbed separatrix,
and introduce the \emph{splitting function}
$\theta\in\T^3\mapsto\M(\theta)\in\R^3$,
providing the vector distance between the whiskers on this section,
along the complementary directions.
In this way, one obtains a measure for the maximal splitting distance
as the maximum of the function~$\abs{\M(\theta)}$.
On the other hand, in suitable coordinates the splitting function
is the gradient of a scalar function called \emph{splitting potential}
\cite{Eliasson94,DelshamsG00},
\[
  \M(\theta)=\nabla\Lc(\theta),
\]
which implies that there always exist homoclinic orbits, which correspond to
the zeros of $\M(\theta)$, i.e.~the critical points of~$\Lc(\theta)$.

In order to provide a first order approximation
to the splitting function, with respect to the parameter $\mu$,
it is very usual to apply the \emph{Poincar\'e--Melnikov method}, introduced by
Poincar\'e in his memoir \cite{Poincare90} and rediscovered much later by
Melnikov and Arnold \cite{Melnikov63,Arnold64}.
This method provides an approximation
\beq\label{eq:melniapproxM}
  \M(\theta)=\mu M(\theta)+\Ord(\mu^2)
\eeq
given by the (vector) \emph{Melnikov function} $M(\theta)$,
defined by an integral (see for instance \cite{Treschev94,DelshamsG00}).
As a result, one obtains asymptotic estimates for the
maximum of the function~$\abs{\M(\theta)}$, provided $\mu$ is small enough.
In fact, the Melnikov function can also be written as the gradient
of a scalar function called the \emph{Melnikov potential}:
$M(\theta)=\nabla L(\theta)$.

However, the case of fast frequencies $\omega_\eps$
as in~(\ref{eq:omega_eps}), with a perturbation of order $\mu=\eps^r$,
for a given $r$ as small as possible,
turns out to be, as said before, a \emph{singular problem}.
The difficulty comes from the fact that the Melnikov function $M(\theta)$
is exponentially small in $\eps$,
and the Poincar\'e--Melnikov method can be directly applied
only if one assumes that $\mu$ is exponentially small with respect to $\eps$
(see for instance \cite{DelshamsG01} for more details).
In order to validate the method in the case $\mu=\eps^r$,
one has to ensure that the error term is also exponentially small,
and that the Poincar\'e--Melnikov approximation dominates it.
To overcome such a~difficulty in the study of the exponentially small splitting,
Lazutkin introduced in \cite{Lazutkin03}
the use of parameterizations of the whiskers on a complex strip
(whose width is defined by the singularities
of the unperturbed parameterized separatrix)
by periodic analytic functions, together with flow-box coordinates.
This tool was initially developed for the Chirikov standard map
\cite{Lazutkin03}, and allowed several authors
to validate the Poincar\'e--Melnikov method
for Hamiltonians with one and a~half degrees of freedom (with only 1~frequency)
\cite{HolmesMS88,Scheurle89,DelshamsS92,DelshamsS97,Gelfreich97}
and for area-preserving maps \cite{DelshamsR98}.

Later, those methods were extended to the case of whiskered tori
with 2~frequencies: $\omega=(1,\Omega)$.
In this case, the arithmetic properties of the frequencies
play an important role in the exponentially small asymptotic estimates of the
splitting function, due to the presence of \emph{small divisors}
of the form $k_1+k_2\Omega$ for integer numbers $k_1$, $k_2$.
Such arithmetic properties can be carefully studied with the help of
the standard theory of \emph{continued fractions}.
The role of the small divisors in the estimates of the splitting
was first noticed by Lochak \cite{Lochak90}
(who obtained an upper bound with an exponent coinciding with Nekhoroshev
resonant normal forms~\cite{Nekhoroshev77}),
and also by Sim\'o \cite{Simo94}
(generalizing an averaging procedure introduced in \cite{Neishtadt84}).
Analogous estimates could also be obtained
from a careful averaging out of the fast angular variables \cite{ProninT00},
at least concerning sharp upper bounds of the splitting.

On the other hand, a numerical detection of asymptotic estimates
was carried out in \cite{Simo94},
and they were rigorously proved in \cite{DelshamsGJS97}
for the quasiperiodically forced pendulum,
assuming a polynomial perturbation in the coordinates associated
to the pendulum. A more general (meromorphic) perturbation
was considered in \cite{GuardiaS12}.
It is worth mentioning that,
in some cases, the Poincar\'e--Melnikov method does not predict correctly
the size of the splitting, as shown in \cite{BaldomaFGS12},
where a Hamilton--Jacobi method is instead used.
This method had previously been used in \cite{Sauzin01,LochakMS03,RudnevW00}.
Similar asymptotic results were obtained in \cite{DelshamsG04}
for the concrete case of the famous \emph{golden ratio} $\Omega=(\sqrt5-1)/2$,
and in \cite{DelshamsGG14c} for the case of
the \emph{silver ratio} $\Omega=\sqrt2-1$,
and generalized in \cite{DelshamsGG16} to any \emph{quadratic} frequency ratio,
and in \cite{DelshamsGG14b} to any frequency ratio \emph{of constant type},
i.e.~with bounded partial quotients.
Very recent results for frequency vectors with unbounded partial quotients
can be found in \cite{FontichSV18a,FontichSV18b}.

In this paper, we consider a 3-dimensional torus with a frequency vector
$\omega$ as in~(\ref{eq:omega_eps})
whose ratios generate a \emph{complex cubic field}
(for short, we say a cubic vector ``of complex type'').
An important difference with respect to the 2-dimensional case
is that in the 3-dimensional case there
is no standard theory of continued fractions
allowing a simple analysis of the small divisors.
As a paradigmatic example, we consider $\omega=(1,\Omega,\Omega^2)$
where $\Omega\approx0.682328$ is the real number
satisfying $\Omega^3=1-\Omega$,
which has been called the \emph{cubic golden number}
(see for instance \cite{HardcastleK00}).
Other famous exemples have been considered in \cite{Chandre02}
(see also \cite{Lochak92} for an account of examples
and results concerning cubic frequencies).

Our goal is to develop a methodology,
based on iteration matrices from a result by Koch \cite{Koch99}
(see Section~\ref{sect:resonantseq1})
allowing us to study the resonances of the given cubic frequency vector.
As a result, we obtain asymptotic estimates for the maximal splitting distance,
whose dependence on $\eps$ is described
by a positive \emph{piecewise-smooth function} denoted~$h_1(\eps)$
(see Theorem~\ref{thm:main}).
In this paper it is proved for the first time
that this function is \emph{quasiperiodic} (and \emph{not periodic})
with respect to $\ln\eps$ with two frequencies $\alpha_1$ and $\alpha_2$,
and its behavior depends \emph{strongly} on the arithmetic properties
of the cubic frequency vector~$\omega$.
In particular, we show that the function $h_1(\eps)$
can be constructed explicitly from the study of the \emph{quasi-resonances}
of the frequency vector $\omega$,
and we can also determine explicitly the frequencies $\alpha_1$ and $\alpha_2$,
as well as upper and lower bounds for $h_1(\eps)$.
In this way, we provide an indication of the complexity of the dependence
on $\eps$ of the splitting.

Such results were partially established in the announcement
\cite{DelshamsGG14a} with a~parallel study of the quadratic and cubic cases
(with~2 and 3~frequencies, respectively),
obtaining also exponentially small estimates
for the maximal splitting distance,
showing the periodicity of the function $h_1(\eps)$
with respect to $\ln\eps$ in the quadratic case
(we also stress that this function becomes a constant in the
case of only 1~frequency, see for instance \cite{DelshamsS97}).
Nevertheless, in \cite{DelshamsGG14a} the quasiperiodicity
of the function $h_1(\eps)$ in the cubic case was only conjectured.

We point out that the aim of this paper is to obtain estimates
for the \emph{maximal splitting distance},
like in our paper~\cite{DelshamsGG14b}
where we considered frequencies of constant type for a 2-dimensional torus.
This is in constrast with most of the papers quoted
in the previous paragraphs, which rather focus their attention on the
\emph{transversality} of the splitting.
The study of the transversality could also be carried out with the methodology
developed here, by means of a more accurate study, as done
in \cite{DelshamsG04,DelshamsGG14c,DelshamsGG16} for the quadratic case
(see remark~\ref{rk:transv} after Theorem~\ref{thm:main}).
We stress that, for some purposes, it is not necessary
to establish the transversality of the splitting,
and it can be enough to provide estimates of the maximal splitting distance.
Indeed, such estimates imply the existence of splitting between the invariant
manifolds, which provides a strong indication of the non-integrability
of the system near the given torus, and opens the door to the application
of topological methods \cite{GideaR03,GideaL06} for the study of
Arnold diffusion in such systems.

\subsection{Setup}

Here we describe the nearly-integrable
Hamiltonian system under consideration.
In particular, we study a \emph{singular} or
\emph{weakly hyperbolic} (\emph{a priori stable})
Hamiltonian with 4 degrees of freedom possessing a
3-dimensional whiskered torus with fast frequencies.
In canonical coordinates
$(x,y, \varphi, I)\in\T\times \R \times \T^3
\times \R^3$, with the symplectic form
$\df x \wedge\df y +\df\varphi\wedge\df I$, the Hamiltonian is defined by
\bea
  \label{eq:HamiltH}
  &&H(x,y, \varphi, I) = H_0 (x,y, I) + \mu H_1(x, \varphi),
\\
  \label{eq:HamiltH0}
  &&H_0 (x, y, I) =
  \langle \omega_\eps, I\rangle + \frac{1}{2} \langle\Lambda I, I\rangle
  +\frac{y^2}{2} + \cos x -1,
\\
  \label{eq:HamiltH1}
  &&H_1 (x, \varphi)= h(x) f(\varphi).
\eea
Our system has two parameters $\eps>0$ and $\mu$, linked
by a relation $\mu=\eps^r$, $r>0$ (the smaller $r$ the better).
Thus, if we consider $\eps$ as the unique parameter, we have
a singular problem for $\eps\to 0$.
See \cite{DelshamsG01} for a discussion about singular and regular problems.

Recall that we are assuming a vector of fast frequencies
$\omega_\eps = \omega/\sqrt{\eps}$ with a cubic vector $\omega\in\R^3$
of \emph{``complex type''}, as introduced in~(\ref{eq:omega_eps}).
It is a well-known property
(and we prove it in Section~\ref{sect:resonantseq2};
see also \cite[\S V.3]{Cassels57} or \cite[\S II.4]{Schmidt80})
that any (complex or totally real) cubic vector
satisfies a \emph{Diophantine condition}
\beq\label{eq:DiophCond}
\abs{\scprod k\omega}\ge\frac{\gamma}{\abs k^2}\,,
\qquad \forall k\in \Z^3\setminus\pp0,
\eeq
with some $\gamma>0$
(the exponent $2$ in this condition is the minimal one
among vectors in $\R^3$).
We also assume in~(\ref{eq:HamiltH0}) that $\Lambda$ is a symmetric
($3\times3$)-matrix, such that $H_0$ satisfies the condition of
\emph{isoenergetic nondegeneracy}
\beq
    \det \left(
    \begin{array}{cc}
    \Lambda & \omega\\
    \omega\tp & 0
    \end{array}
    \right) \neq 0.
\label{eq:isoenerg}
\eeq

For the perturbation $H_1$ in~(\ref{eq:HamiltH1}),
we deal with the following analytic periodic functions,
\beq\label{eq:hf}
  h(x) = \cos x,
  \quad
  f(\varphi)
  =\sum_{k\in\Zc}f_k\cos(\langle k, \varphi\rangle-\sigma_k),
  \quad \textrm{with} \quad f_k=\ee^{-\rho\abs k}
  \ \ \textrm{and} \ \ \sigma_k\in\T,
\eeq
where we introduce, in order to avoid repetitions in the Fourier series,
the set
\beq\label{eq:calZ}
  \Zc=\{k\in\Z^3:
  k_2\ge1\ \textrm{or}\ (k_2=0,k_3\ge1)\ \textrm{or}\ (k_2=k_3=0,k_1\ge0)\},
\eeq
with $k=(k_1,k_2,k_3)$
(the specific choice of $k_2$ being positive, which is not relevant,
allows us to agree with the definition of the set $\Pc$ in~(\ref{eq:calP})).
Notice that, for any couple $\pm k$ of integer vectors,
only one of them belongs to $\Zc$.
The constant $\rho>0$ gives the complex width of analyticity
of the function $f(\varphi)$.
Concerning the phases $\sigma_k$, they can be chosen arbitrarily
for the purpose of this paper.

To justify the form of the perturbation $H_1$
chosen in~(\ref{eq:HamiltH1}) and~(\ref{eq:hf}),
we stress that it makes easier the explicit
computation of the Melnikov potential, which
is necessary in order to show that it dominates
the error term in~(\ref{eq:melniapproxM}),
and therefore to establish the existence of splitting.
Moreover, the assumption that all coefficients $f_k$
in the Fourier expansion~(\ref{eq:hf}) with respect to $\varphi$
are nonzero and have an exponential decay,
is usual in the literature
(see for instance \cite{FontichSV18a,FontichSV18b}), and ensures
that the study of the dominant harmonics of the Melnikov potential
can be carried out directly from the arithmetic properties
of the frequency vector~$\omega$.
Indeed, such dominant harmonics
correspond to the integer vectors $k$ providing an approximate equality
in~(\ref{eq:DiophCond}), i.e.~giving the ``smallest'' divisors
(relatively to the size of $\abs k$).
We call \emph{primary resonances} of $\omega$ to such vectors~$k$,
and \emph{secondary resonances} to the rest of quasi-resonances
(see Section~\ref{sect:cubicfreq} for details).
In this way, the choice of the coefficients $f_k$ in~(\ref{eq:hf})
allows us to emphasize the dependence of the splitting
on the arithmetic properties of $\omega$.

It is worth remarking that, once we know the primary resonances
for the given frequency vector $\omega$, we do not need
all the coefficients~$f_k$ to be different from zero in~(\ref{eq:hf}),
but only the ones corresponding to primary resonances.
On the other hand, since our method is completely constructive, other choices
of concrete harmonics $f_k$ could also be considered
(like $f_k=\abs k^m\ee^{-\rho\abs k}$), simply at the cost of more cumbersome
computations in order to determine the dominant harmonics
of the Melnikov potential.

We also remind that the Hamiltonian
defined in~(\ref{eq:HamiltH}--\ref{eq:hf}) is paradigmatic,
since it is a generalization of the famous Arnold's example
(introduced in \cite{Arnold64} to illustrate the transition
chain mechanism in Arnold diffusion).
It provides a model for the behavior of
a nearly-integrable Hamiltonian system near a single resonance
(see \cite{DelshamsG01} for a motivation),
and has often been considered in the literature
(see for instance \cite{GallavottiGM99b,ProninT00,LochakMS03,DelshamsGS04}).

Let us describe the invariant tori and whiskers,
as well as the splitting and Melnikov functions.
First, it is clear that the unperturbed system given by $H_0$
(that corresponds to $\mu=0$) is separable,
and consists of the pendulum given by $P(x,y)=y^2/2+\cos x-1$, and
3~rotors with fast frequencies:
\ $\dot{\varphi}= \omega_\eps +\Lambda I$, \ $\dot I=0$.
\ The pendulum has a hyperbolic equilibrium at the origin, with
separatrices that correspond to the curves given by $P(x,y)= 0$.
We parameterize the upper separatrix of the pendulum as
$(x_0(s), y_0(s))$, $s\in \R$, where
\[
  x_0(s)=4\arctan\ee^s, \qquad y_0(s) = \frac{2}{\cosh s}\,.
\]
Then, the lower separatrix has the parametrization $(x_0(-s),-y_0(-s))$.
For the rotors system $(\varphi, I)$, the solutions are
$I= I_0$, $\varphi = \varphi_0+t(\omega_\eps+\Lambda I_0)$.
Consequently, the Hamiltonian $H_0$ has a 3-parameter family of 3-dimensional
whiskered tori: in coordinates $(x,y,\varphi,I)$,
each torus can be parameterized as
\[
  \Tc_{I_0}:\qquad(0,0,\theta,I_0),\quad\theta\in\T^3,
\]
and the inner dynamics on each torus is
$\dot\theta = \omega_\eps+\Lambda I_0$.
Each invariant torus has a \emph{homoclinic whisker},
i.e.~coincident 4-dimensional stable and unstable invariant manifolds,
which can be parameterized as
\beq\label{eq:W0}
  \W_{I_0}:\qquad(x_0(s), y_0(s),\theta,I_0),
  \quad s\in \R, \ \theta\in\T^3,
\eeq
with the inner dynamics given by
$\dot{s}=1$, $\dot{\theta}=\omega_\eps+\Lambda I_0$.

In fact, the collection of the whiskered tori
for all values of $I_0$ is
a 6-dimensional \emph{normally hyperbolic invariant manifold},
parameterized by $(\theta,I)\in\T^3\times\R^3$.
This manifold has a 7-dimensional homoclinic manifold,
which can be parameterized by $(s,\theta,I)$,
with inner dynamics
$\dot s=1$, $\dot\theta=\omega_\eps+\Lambda I$, \ $\dot I=0$.
We stress that this approach is usually considered
in the study of Arnold diffusion
(see for instance \cite{DelshamsLS06}).

Among the family of whiskered tori and homoclinic whiskers,
we are going to focus our attention on the torus $\Tc_0$,
whose frequency vector is $\omega_\eps$ as in~(\ref{eq:omega_eps}),
and its associated homoclinic whisker $\W_0$.

When adding the perturbation $\mu H_1$, for $\mu\ne0$ small enough
the \emph{hyperbolic KAM theorem} can be applied
thanks to the Diophantine condition~(\ref{eq:DiophCond})
and to the isoenergetic nondegeneracy~(\ref{eq:isoenerg}).
For $\mu$ small enough, the whiskered torus persists
with some shift and deformation, as a perturbed torus $\Tc=\Tc^{(\mu)}$,
as well as its local whiskers $\W_\loc=\W^{(\mu)}_\loc$
(precise statements can be found, for instance,
in \cite{Niederman00,DelshamsGS04}).

The local whiskers can be extended along the flow,
but in general for $\mu\ne0$ the \emph{(global) whiskers}
do not coincide anymore, and one expects the existence of splitting between
the (4-dimensional) stable and unstable whiskers,
denoted $\W^\st=\W^{\st,(\mu)}$ and $\W^\ut=\W^{\ut,(\mu)}$ respectively.
Using \emph{flow-box coordinates}
(see \cite{DelshamsGS04}, where the $n$-dimensional case is considered)
in a neighbourhood containing a piece of
both whiskers (away from the invariant torus),
one can introduce parameterizations of the perturbed whiskers,
with parameters $(s,\theta)$ inherited
from the unperturbed whisker~(\ref{eq:W0}),
and the inner dynamics
\[
  \dot s=1, \qquad \dot{\theta}=\omega_\eps.
\]
Then, the distance between the stable whisker $\W^\st$ and
the unstable whisker $\W^\ut$ can be measured by comparing
such parameterizations along the complementary directions.
The number of such directions is~4 but,
due to the energy conservation, it is enough to consider 3~directions,
say the ones related to the action coordinates~$I$.
In this way, one can introduce a (vector) \emph{splitting function},
with values in $\R^3$, as the difference of the parameterizations
$\J^{\st,\ut}(s,\theta)$ of (the action components of)
the perturbed whiskers $\W^\st$ and $\W^\ut$.
Initially this splitting function depends on $(s,\theta)$,
but it can be restricted to a transverse section
by considering a fixed $s$, say $s=0$,
and we can define as in \cite[\S5.2]{DelshamsG00} the splitting function
\beq\label{eq:defM}
  \M(\theta):=\J^\ut(0,\theta)-\J^\st(0,\theta),
  \quad\theta\in\T^3.
\eeq

Applying the Poincar\'e--Melnikov method,
the first order approximation~(\ref{eq:melniapproxM}) of the
splitting function is given by
the (vector) \emph{Melnikov function} $M(\theta)$,
which is the gradient of the (scalar) \emph{Melnikov potential}:
\ $M(\theta)=\nabla L(\theta)$.
\ The latter one can be defined as an integral:
we consider any homoclinic trajectory of
the unperturbed homoclinic whisker~$\W_0$ in~(\ref{eq:W0}),
starting on the section $s=0$, and the trajectory on the torus $\Tc_0$
to which it is asymptotic as $t\to\pm\infty$,
and we substract the values of the perturbation~$H_1$ on the two trajectories.
This gives an absolutely convergent integral,
which depends on the initial phase~$\theta\in\T^3$
of the considered trajectories:
\bea
  \nonumber
  L(\theta)
  &:=
  &-\int_{-\infty}^{\infty}
    [H_1(x_0(t),\theta+t\omega_\eps)-H_1(0,\theta+t\omega_\eps)]\,\df t
\\
  \label{eq:L}
  &=
  &-\int_{-\infty}^{\infty}
    [h(x_0(t))-h(0)] f(\theta+t\omega_\eps)\,\df t,
\eea
where we have taken into account the specific form~(\ref{eq:HamiltH1})
of the perturbation.

Our choice of the pendulum $P(x,y)=y^2/2+\cos x-1$ in~(\ref{eq:HamiltH0}),
whose separatrix has simple poles,
makes it possible to use the method of residues in order to compute
the coefficients $L_k$ of the Fourier expansion of
the Melnikov potential~$L(\theta)$,
and hence for the coefficients of the Melnikov function:
$\abs{M_k}=\abs kL_k$.
Such coefficients turn out to be exponentially small in $\eps$
(see their expression in Section~\ref{sect:gn}).
For each value of~$\eps$ only the \emph{dominant harmonic},
corresponding to some index $k=S_1(\eps)$,
is relevant in order to provide asymptotic estimates
for the maximum value of the Melnikov function
(of course, a few dominant harmonics may have to be considered
near some transition values of $\eps$,
at which changes in the dominance take place).
Due to the exponential decay of the Fourier coefficients of $f(\varphi)$
in~(\ref{eq:hf}),
it is not hard to study such a dominance and its dependence on $\eps$.

In order to give asymptotic estimates for the maximal splitting distance,
the estimates obtained for the Melnikov function $M(\theta)$
have to be validated also for the splitting function~$\M(\theta)$.
The difficulty in the application of
the Poincar\'e--Melnikov approximation~(\ref{eq:melniapproxM}),
due to the exponential smallness in $\eps$
of the function $M(\theta)$ in our singular case $\mu=\eps^r$,
can be solved by obtaining upper bounds (on a complex domain)
for the \emph{error term} in~(\ref{eq:melniapproxM}),
showing that, if $r>r^*$ with a suitable $r^*$, its Fourier coefficients
are dominated by the coefficients of $M(\theta)$
(see also~\cite{DelshamsGS04}).

We stress that our approach can also be directly applied to other
classical 1-degree-of-freedom Hamiltonians $P(x,y)=y^2/2+V(x)$,
with a potential $V(x)$ having a unique nondegenerate maximum,
although the use of residues becomes more cumbersome when the
complex parameterization of the separatrix has poles of higher orders
(see some examples in \cite{DelshamsS97}).

\subsection{Main result}

For the Hamiltonian system~\mbox{(\ref{eq:HamiltH}--\ref{eq:hf})}
with the 2 parameters linked by
$\mu=\eps^r$, $r>r^*$ (with some suitable $r^*$),
and a cubic frequency vector of complex type $\omega$
as in~(\ref{eq:omega_eps}),
our main result provides an exponentially small \emph{asymptotic estimate}
for the \emph{maximal distance} of splitting,
given in terms of the maximum size in modulus of the
splitting function $\M(\theta)$,
and this estimate is valid for all $\eps$ sufficiently small.

With our approach, the Poincar\'e--Melnikov method can be validated
for an exponent $r>r^*$ with $r^*=3$,
although a lower value of $r^*$ can be given in some particular cases
(see remark~\ref{rk:hf2} after Theorem~\ref{thm:main}).
However, such values of $r^*$ are not optimal and could be improved
using other methods, like the parametrization of the
whiskers as solutions of Hamilton--Jacobi equation
(see for instance \cite{LochakMS03,BaldomaFGS12}).
In this paper, the emphasis is put on the extension of the methods
and results from the 2-dimensional quadratic case
to the 3-dimensional cubic case,
rather than on the improvement of the value of $r^*$.

Due to the form of $f(\varphi)$ in~(\ref{eq:hf}),
the Melnikov potential $L(\theta)$
is readily presented in its Fourier series (see Section~\ref{sect:gn}),
with coefficients $L_k=L_k(\eps)$
which are exponentially small in $\eps$.
We use this expansion of $L(\theta)$ in order to detect its
\emph{dominant harmonic} $k=S_1(\eps)$ for every given $\eps$.
Such a dominance is also valid for the Melnikov function
$M(\theta)$, since the size of their Fourier coefficients
$M_k$ (vector) and $L_k$ (scalar) is directly related:
$\abs{M_k}=\abs k\,L_k$, $k\in\Zc$
(recall the definition of $\Zc$ in~(\ref{eq:calZ})).

As shown in Section~\ref{sect:technical}, in order to obtain an
asymptotic estimate for the maximum value of $\M(\theta)$,
i.e.~for the distance of splitting, for most values of $\eps$
it is enough to consider the (unique) first dominant harmonic $S_1(\eps)$
of the Melnikov function $M(\theta)$,
whose size behaves like $\exp\{-h_1(\eps)/\eps^{1/6}\}$,
being described by a (positive) function $h_1(\eps)$
that is carefully studied in this paper.
To ensure that the dominant harmonic of $M(\theta)$ corresponds
to the dominant harmonic of the splitting function $\M(\theta)$,
one has to carry out an accurate control
of the error term in~(\ref{eq:melniapproxM}).
In this way, using estimates for the size
of the dominant harmonic, as well as for all the remaining harmonics,
one can prove that the dominant harmonic
is large enough and provides an approximation
to the maximum size of the whole splitting function
(see also \cite{DelshamsGG14a,DelshamsGG14b,DelshamsGG16}).

\begin{figure}[!b]
  \centering
  \subfigure{
    \includegraphics[width=0.24\textwidth,angle=-90]{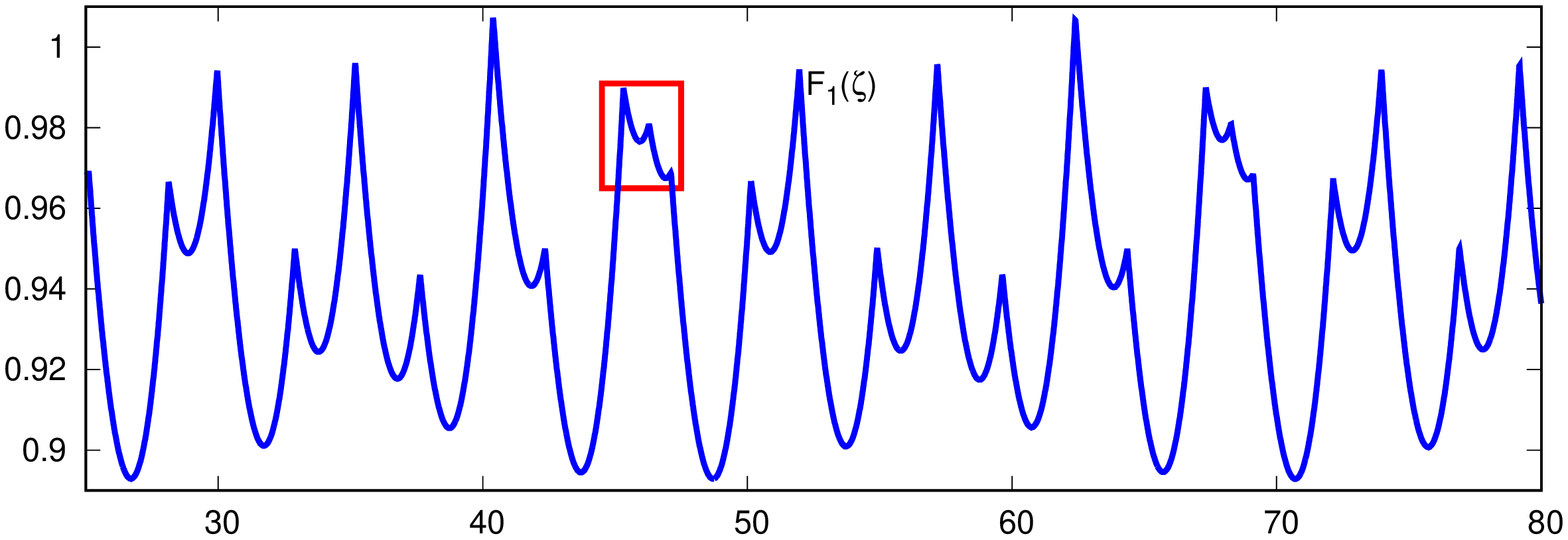}}
  \subfigure{
    \includegraphics[width=0.23\textwidth,angle=-90]{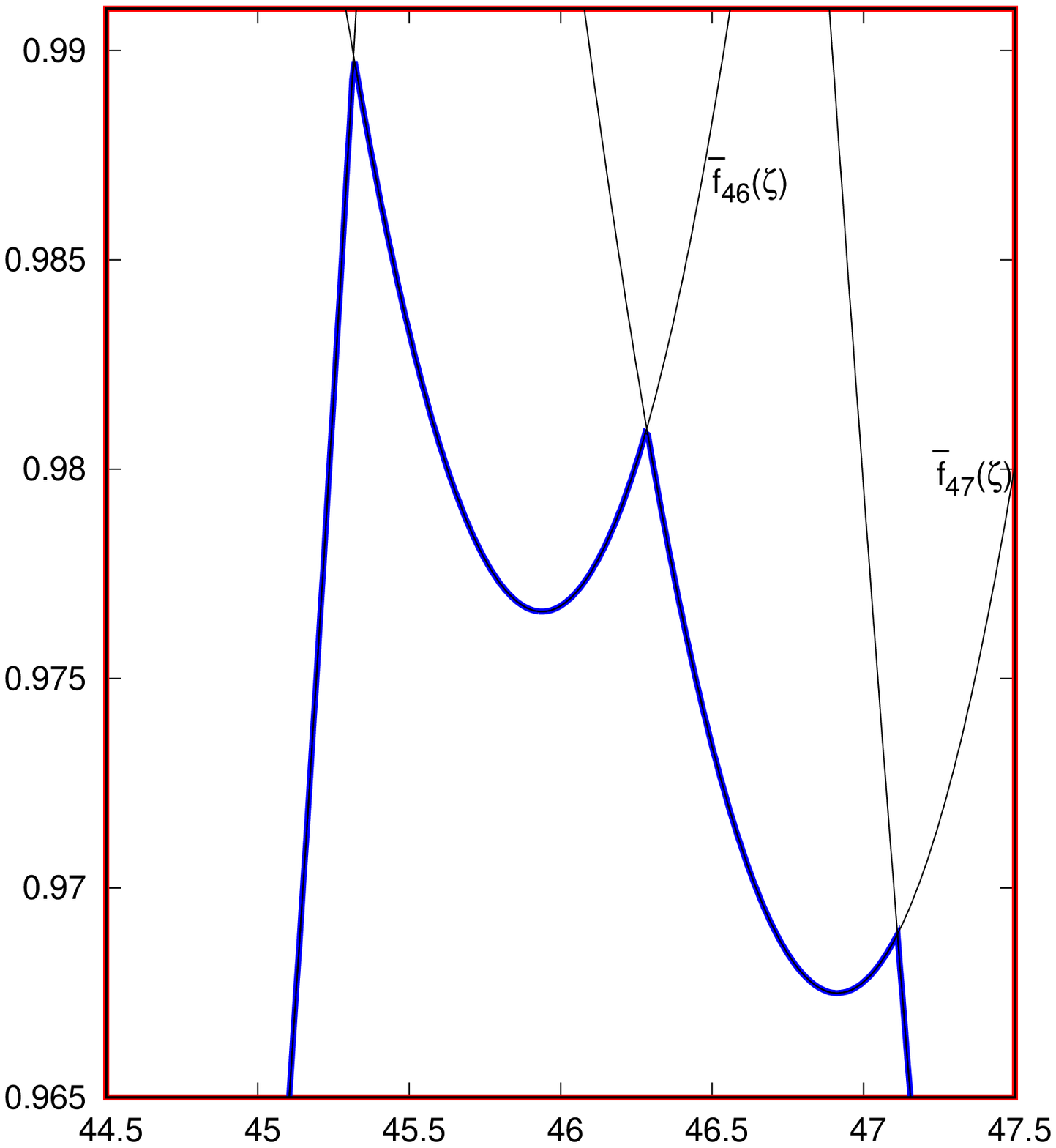}}
  \caption{\small\emph{
    Graph of the function $h_1(\eps)=F_1(\zeta)$ in the exponent
    of~(\ref{eq:main}), for the cubic golden vector
    (see Section~\ref{sect:cubicgolden}),
    in the logarithmic variable $\zeta\sim\ln(1/\eps)$
    (see~(\ref{eq:logchange}) for a precise definition),
    as the minimum of the functions $\bbf_n(\zeta)$
    (see Sections~\ref{sect:h1} and~\ref{sect:h1_cubicgolden}).}}
\label{fig:h1}
\end{figure}

However, one has to consider at least two harmonics
for $\eps$ near to some \emph{``transition values''},
at which a change in the dominant harmonic occurs and,
consequently, two (or more) harmonics having similar sizes
can be considered as the dominant ones.
In this case, the size of the splitting function
can also be determined from the dominant harmonics,
although such transition values turn out to be \emph{corners}
of the function $h_1(\eps)$
(see the theorem below, and Figure~\ref{fig:h1}).

The determination of the dominant harmonics, and
hence the dependence on $\eps$ of the size of the splitting
and the function $h_1(\eps)$,
are closely related to the arithmetic properties of the
frequency vector $\omega$ in~(\ref{eq:omega_eps}), since
the integer vectors $k\in\Zc$ associated to the dominant harmonics
can be found, for any $\eps$, among
the main quasi-resonances of $\omega$,
i.e.~the vectors $k$ giving the
``smallest'' divisors $\abs{\scprod k\omega}$
(relatively to the size of $\abs k$).
In Section~\ref{sect:cubicfreq}, we develop
a methodology for a complete study of the \emph{resonant properties}
of cubic frequency vectors (of complex type),
which is one of the main goals of this paper.
This methodology relies on the classification
of the integer vectors~$k$ into \emph{``resonant sequences''}
(see Section~\ref{sect:resonantseq1} for definitions).
Among them, the sequence of \emph{primary resonances}
corresponds to the vectors $k$ which fit best
the Diophantine condition~(\ref{eq:DiophCond}),
and the vectors $k$ belonging to the remaining sequences
are called \emph{secondary resonances}.
In this way, we can also determine
the (positive) \emph{asymptotic Diophantine constant},
\beq\label{eq:gammam}
  \gamma^-:=\liminf_{\abs k\to\infty}\abs{\scprod k\omega}\cdot\abs k^2.
\eeq
This approach, already announced in \cite{DelshamsGG14a}
for 3-dimensional cubic frequency vectors,
generalizes the one introduced in \cite{DelshamsG03}
for 2-dimensional quadratic frequency vectors.

For most values of $\eps$, the dominant harmonic is given by
an integer vector~$k$ associated to a primary resonance,
but for some intervals of $\eps$ the secondary resonances
may have to be taken into account giving rise to
a more involved function $h_1(\eps)$.
Nevertheless, for some cubic frequency vectors $\omega$ in~(\ref{eq:omega_eps})
such as the \emph{cubic golden vector},
the function~$h_1(\eps)$ can be defined using only the primary resonances
(see Sections~\ref{sect:cubicgolden} and~\ref{sect:h1_cubicgolden}).

In order to generate the resonant sequences,
we use a result by Koch \cite{Koch99}, ensuring the existence of
a \emph{unimodular} ($3\times3$)-matrix $T$
(i.e.~with integer entries and determinant $\pm1$),
having $\omega$ as an eigenvector with the associated eigenvalue
\beq\label{eq:lambda0}
  \lambda>1.
\eeq
Altough there exist an infinity of matrices $T$ fitting Koch's result,
we establish in Section~\ref{sect:resonantseq1} a canonical choice for it
(see Proposition~\ref{prop:uniqueT}), and we write it as $T=T(\omega)$.

The eigenvalue $\lambda=\lambda(\omega)$ is also a cubic irrational number
and belongs to $\Q(\Omega)$. Hence it also has complex conjugates,
which can be written in the form
\beq\label{eq:phi0}
  \lambda_2,\ol\lambda_2=\frac1{\sqrt\lambda}\,\ee^{\pm\ii\pi\cdot\phi},
  \qquad
  0<\phi<1,
\eeq
and $\phi=\phi(\omega)$ is an irrational number
(see Section~\ref{sect:resonantseq1}).

For a concrete cubic frequency vector $\omega$,
it is not too hard to find the Koch's matrix $T=T(\omega)$
(see Section~\ref{sect:resonantseq1} for a procedure,
and Section~\ref{sect:cubicgolden} for its application
to the concrete case of the cubic golden vector).
We point out that, for the quadratic 2-dimensional case $\omega=(1,\Omega)$,
a systematic algorithm providing an analogous ($2\times2$)-matrix $T$
was developed in \cite{DelshamsGG16},
from the continued fraction of the frequency ratio $\Omega$
(which is eventually periodic for quadratic numbers).
An extension of this algorithm to the cubic case would require a further study
(possibly using some of the existing multidimensional continued fraction
theories), and is not carried out here.

Assuming that the matrix $T$ is known,
the key point is that the iteration of the matrix $U=(T^{-1})\tp$
from an initial (\emph{``primitive''}) vector allows us
to generate any resonant sequence (see the definition~(\ref{eq:sqn})).
In this way, we can construct the resonant sequences
allowing us to detect the dominant harmonics of the Melnikov potential
and, consequently, asymptotic estimates for the maximal splitting distance.

Next, we establish the \emph{main result} of this work,
which generalizes to the complex cubic case
the results obtained in \cite{DelshamsG04,DelshamsGG16} for the quadratic case.
The result given below provides exponentially small asymptotic estimates
for the maximal distance of splitting, as $\eps\to0$,
given by the maximum of $\abs{\M(\theta)}$, $\theta\in\T^3$.
In such asymptotic estimates, the dependence on $\eps$
is mainly described by the exponent $1/6$, and by the function $h_1(\eps)$.
This is a positive function, \emph{quasiperiodic}
with respect to $\ln\eps$ and \emph{piecewise-smooth}
and, consequently, it has a finite number of \emph{corners}
(i.e.~jump discontinuities of the derivative) in any given interval.
As we can see from the statement of the theorem,
the numbers $\lambda$ and $\phi$
introduced in~(\ref{eq:lambda0}--\ref{eq:phi0})
play an essential role in the quasiperiodicity of the function $h_1(\eps)$,
since they provide directly the two frequencies
$3\ln\lambda$ and $3\ln\lambda\cdot\phi$,
and the fact that $\phi$ is irrational ensures that the function $h_1(\eps)$
is \emph{not periodic}, which makes a difference with respect to the
quadratic case considered in \cite{DelshamsGG16}.

For any given cubic vector $\omega$ (of complex type),
the function $h_1(\eps)$ can be explicitly constructed
(see Section~\ref{sect:h1}).
However, its (piecewise) expression can be very complicated.
Its graph is shown in Figure~\ref{fig:h1}
(where a \emph{logarithmic scale} for $\eps$ is used),
for the concrete case of the cubic golden frequency vector.
The oscillatory behavior of the function $h_1(\eps)$
depends strongly on the arithmetic properties of $\omega$.

For positive quantities, we use the notation $f\sim g$
if we can bound $c_1g\le f\le c_2g$ with constants
$c_1,c_2>0$ not~depending on $\eps$, $\mu$.

\begin{theorem}[\emph{main result}]\label{thm:main}
Assume the conditions described for the
Hamiltonian~\mbox{(\ref{eq:HamiltH}--\ref{eq:hf})},
with a cubic frequency vector $\omega=(1,\Omega,\wtl\Omega)$ of complex type
as in~(\ref{eq:omega_eps}),
that $\eps$ is small enough and that $\mu=\eps^r$, $r>3$.
Then, for the splitting function $\M(\theta)$ we have:
\beq\label{eq:main}
  \max_{\theta\in\T^3}\abs{\M(\theta)}
  \sim \frac{\mu}{\eps^{1/3}}
  \exp \pp{- \frac{C_0 h_1 (\eps)}{\eps^{1/6}}}.
\eeq
The function $h_1(\eps)$, defined in~(\ref{eq:h1}),
is positive, piecewise-smooth, piecewise-convex and quasiperiodic in $\ln\eps$,
with two frequencies $3\ln\lambda$ and $3\ln\lambda\cdot\phi$,
where $\lambda=\lambda(\omega)$ and $\phi=\phi(\omega)$ are the numbers
introduced in~(\ref{eq:lambda0}--\ref{eq:phi0}).
It satisfies for $\eps>0$ lower and upper bounds $J^-_0\le h_1(\eps)\le J^+_1$,
where the values $J^-_0=J^-_0(\omega)$ and $J^+_1=J^+_1(\omega)$
are defined in~(\ref{eq:J01}).
On the other hand, $C_0=C_0(\omega,\rho)$ is a positive constant
defined in~(\ref{eq:beta_gk}).
\end{theorem}

\bremarks
\item\label{rk:supremum}
As a consequence of this theorem, replacing $h_1(\eps)$ by
its supremum value $J^*_1$ ($\le J^+_1$, see also Section~\ref{sect:qp}),
we get the following \emph{sharp lower bound} for
the maximal splitting distance:
\[
  \max_{\theta\in\T^3} |\M(\theta)|
  \ge \frac{c\mu}{\eps^{1/3}}
  \exp \pp{-\frac{C_0J^*_1}{\eps^{1/6}}},
\]
where $c$ is a constant.
This may be enough, if our aim is only to prove the existence of splitting
of separatrices, without giving an accurate description for it.
\item\label{rk:transv}
Our approach can also be applied to show the existence
of \emph{transverse homoclinic orbits},
associated to simple zeros~$\theta^*$ of the splitting function $\M(\theta)$
(or, equivalently, nondegenerate critical points of the splitting potential),
providing an asymptotic estimate for the \emph{transversality} of the
homoclinic orbits, measured by the minimum eigenvalue (in modulus)
of the matrix $\Df\M(\theta^*)$ at each zero of $\M(\theta)$.
Such an asymptotic estimate is exponentially small in $\eps$
as in~(\ref{eq:main}), but the function $h_1(\eps)$
has to be replaced by a greater function $h_3(\eps)$,
also piecewise-smooth and quasiperiodic in $\ln\eps$.
In order to define $h_3(\eps)$, one has to consider
the three most dominant harmonics whose indices
$S_1(\eps),S_2(\eps),S_3(\eps)\in\Zc$ are linearly independent
(this is necessary in order to prove that the zeros $\theta^*$ are simple).
This result on transversality would be valid for \emph{``almost all''} $\eps$
sufficiently small, since one has to exclude a small neighborhood of
some values where the third and the fourth dominant harmonics
have similar sizes, and homoclinic bifurcations could take place.
See \cite{DelshamsGG16} for the analogous situation in the quadratic case,
where only the two most dominant harmonics are necessary.
\item\label{rk:hf2}
The results of Theorem~\ref{thm:main} can be improved under some
particular situations. For instance, if the function $h(x)$ in~(\ref{eq:hf})
is replaced by~$h(x)=\cos x-1$,
then the estimates are valid for $\mu=\eps^r$ with $r>2$ (instead of $r>3$).
The details of this improvement are not given here,
since they work exactly as in~\cite{DelshamsG04}.
\eremarks

\paragr{Organization of the paper}
We start in Section~\ref{sect:cubicfreq} with studying
the arithmetic properties of cubic frequency vectors
$\omega=(1,\Omega,\wtl\Omega)$ (of complex type),
and constructing the iteration matrix $T$.
Next, in Section~\ref{sect:asympt_est} we find
an asymptotic estimate for the dominant harmonic
of the splitting potential,
which allows us to define the function $h_1(\eps)$
and study their general properties.
In order to illustrate our methods,
concrete results for the cubic golden vector are obtained in
Sections~\ref{sect:cubicgolden} (aritmetic properties)
and~\ref{sect:h1_cubicgolden} (the function $h_1(\eps)$).
Finally, in Section~\ref{sect:technical} we provide rigorous bounds
of the remaining harmonics allowing us
to obtain asymptotic estimates for the maximal splitting distance,
as established in Theorem~\ref{thm:main}.

\section{Arithmetic properties of cubic frequencies}
\label{sect:cubicfreq}

\subsection{Iteration matrix for a cubic frequency vector}
\label{sect:resonantseq1}

We consider a \emph{cubic frequency vector} $\omega\in\R^3$,
i.e.~the frequency ratios $\omega_2/\omega_1$ and $\omega_3/\omega_1$
generate a cubic field (an~algebraic number field of degree~3 over $\Q$,
i.e.~its dimension as a vector space over $\Q$ is~3).
In order to simplify our exposition, we assume that $\omega_1=1$,
and hence the vector has the form
\beq\label{eq:omega}
  \omega=(1,\Omega,\wtl\Omega),
\eeq
where $\Omega$ is a \emph{cubic irrational number},
i.e.~its minimum polynomial
(the monic polynomial of minimal degree having $\Omega$ as a root)
has degree~3, and $\wtl\Omega$ belongs to the field $\Q(\Omega)$:
\bea
  \label{eq:rootOmega}
  &&\Omega^3=r_0+r_1\Omega+r_2\Omega^2,
\\
  \label{eq:tlOmega}
  &&\wtl\Omega=a_0+a_1\Omega+a_2\Omega^2,
  \qquad \textrm{with} \quad a_2\ne0,
\eea
where the coefficients $r_j$, $a_j$ are rational.
The number $\wtl\Omega$ is also cubic irrational
(in fact, any number belonging to~$\Q(\Omega)$
is either rational or cubic irrational).
We restrict ourselves to the \emph{complex case}
(also called the \emph{non-totally real case}):
the two conjugates of $\Omega$, as a root
of the polynomial equation~(\ref{eq:rootOmega}), are complex.
This condition can be expressed in terms
of having \emph{negative discriminant},
\[
  \Delta=4r_1^{\,3}+r_1^{\,2}r_2^{\,2}-27r_0^{\,2}-18r_0r_1r_2-4r_0r_2^{\,3}
  \ <\,0.
\]

We denote the conjugates of $\Omega$ as
\beq\label{eq:sigma23a}
  \Omega_2:=\sigma(\Omega)=\sigma_2+\ii\sigma_3\,,
  \quad
  \ol\Omega_2=\bar\sigma(\Omega)=\sigma_2-\ii\sigma_3
\eeq
and, from the standard equalities
\[
  r_2=\Omega+\Omega_2+\ol\Omega_2=\Omega+2\sigma_2,
  \quad
  r_1=-(\Omega\Omega_2+\Omega\ol\Omega_2+\Omega_2\ol\Omega_2)
  =-(2\Omega\sigma_2+\sigma_2^{\,2}+\sigma_3^{\,2})
\]
we see that
\beq\label{eq:sigma23b}
  \sigma_2=\frac12(r_2-\Omega),
  \qquad
  \sigma_3=\frac s2\,\sqrt{-(4r_1+r_2^{\,2})-2r_2\Omega+3\Omega^2}\,,
\eeq
with a concrete sign $s=\pm1$ for $\sigma_3$, that will be chosen later
for convenience (see~(\ref{eq:phi})).

It is clear from~(\ref{eq:tlOmega}) that our cubic frequency vector $\omega$
can be related to the more particular case
\beq\label{eq:omega0}
  \omega^{(0)}=(1,\Omega,\Omega^2)
\eeq
through a linear change:
$\omega=A\,\omega^{(0)}$, with the following matrix belonging
to the \emph{general linear group} $\GL(3,\Q)$,
\beq\label{eq:A}
  A:=\p{\begin{array}{ccc}
          1   &0   &0\\
          0   &1   &0\\
          a_0 &a_1 &a_2
        \end{array}}
\eeq
(for instance, the cubic golden frequency vector considered
in Section~\ref{sect:cubicgolden} has the form~(\ref{eq:omega0})).

It is well-known from algebraic number theory
(see for instance \cite[ch.~II]{StewartT87}
or \cite[ch.~V--VI]{Lang02} as general references)
that there exist unique field isomorphisms
$\sigma:\Q(\Omega)\longrightarrow\Q(\Omega_2)$
and $\bar\sigma:\Q(\Omega)\longrightarrow\Q(\ol\Omega_2)$
such that $\sigma(\Omega)=\Omega_2$ and $\bar\sigma(\Omega)=\ol\Omega_2$.
It is clear that $\sigma$ and $\bar\sigma$ are related
by the ordinary complex conjugacy.
Then, the numbers $\sigma(\wtl\Omega)$ and $\bar\sigma(\wtl\Omega)$
turn to be the conjugates of $\wtl\Omega$,
and they are also complex
(indeed, if they were real, they would coincide and
$\wtl\Omega$ would not be a cubic irrational).

Any cubic frequency vector $\omega\in\R^3$
satisfies a \emph{Diophantine condition},
with the minimal exponent
(see for instance \cite[\S V.3]{Cassels57} or \cite[\S II.4]{Schmidt80}):
\beq\label{eq:diophantine}
  \abs{\scprod k\omega}\ge\frac\gamma{\abs k^2}\,,
  \quad
  \forall k\in\Z^3\setminus\pp0.
\eeq
With this in mind, we define the \emph{``numerators''}
\beq\label{eq:numerators}
  \gamma_k:=\abs{\scprod k\omega}\cdot\abs k^2,
  \qquad
  k\in\Z^3\setminus\pp0,
\eeq
where we use the Euclidean norm: $\abs\cdot=\absd\cdot$
(this allows us to use the properties of the scalar product).
The numerators have $\gamma>0$ as a lower bound.
Our goal is to provide a classification of the integer vectors $k$,
according to the size of $\gamma_k$, in order to find
the \emph{primary resonances}
(i.e.~the integer vectors $k$ for which $\gamma_k$ is smallest, and
hence best fitting the Diophantine condition~(\ref{eq:diophantine})),
and study their separation with respect to the remaining vectors $k$
(i.e.~the secondary resonances).

The key point will be to use the following result by Koch \cite{Koch99}:
for a vector $\omega\in\R^\ell$ whose frequency ratios generate
an algebraic field of degree~$\ell$,
there exists a \emph{unimodular} $(\ell\times\ell)$-matrix $T$
(a square matrix with integer entries and determinant~$\pm 1$)
having $\omega$ as an eigenvector with associated eigenvalue $\lambda$
of modulus~$>1$, and such that the other $\ell-1$ eigenvalues are all simple
and of modulus~$<1$.
This result is valid for any dimension~$\ell$,
and is usually applied in the context of renormalization theory
(see for instance \cite{Koch99,Lopesd02a}), since the iteration of the
matrix $T$ provides successive rational approximations to the direction
of the vector~$\omega$.

For any given cubic frequency vector $\omega$ as in~(\ref{eq:omega}),
we say that a $(3\times3)$-matrix~$T$
is a ``\emph{Koch's matrix for $\omega$}\,''
if~it satisfies the requirements of Koch's result \cite{Koch99}.
It is not hard to find a Koch's matrix for any concrete cubic vector~$\omega$
(see below for a general procedure, and Section~\ref{sect:cubicgolden}
for its application to the concrete case of the cubic golden vector).
It is clear that a Koch's matrix $T$ is not unique, since
any power $\pm T^n$ is also a Koch's matrix.

We can assume that the determinant of $T$ is positive, $\det T=1$,
i.e.~$T$ belonging to the \emph{special linear group} $\SL(3,\Z)$
(otherwise, we can replace $T$ by~$-T$).
For the eigenvalue $\lambda$ associated to the eigenvector~$\omega$,
it is clear that it is real and can be writen as
\beq\label{eq:lambda}
  \lambda=\scprod{T_{(1)}}\omega
  =T_{11}+T_{12}\Omega+T_{13}\wtl\Omega\ \in\ \Q(\Omega)
\eeq
where we denote $T_{(1)}:=(T_{11},T_{12},T_{13})$
(the first row of $T$, considered here as a column vector).
We also see that $\lambda$ is cubic irrational
(otherwise, it would be rational and
the frequency ratios of $\omega$ would also be rational).
The other two eigenvalues of $T$, which are the conjugates of $\lambda$,
are complex (see the argument given above for $\wtl\Omega$),
which implies that $\lambda$ is positive: $\lambda>1$.
We write the conjugates of $\lambda$ in terms of real and imaginary parts:
\beq\label{eq:lambda2}
  \lambda_2:=\sigma(\lambda)=\mu_2+\ii\mu_3\,,
  \quad
  \ol\lambda_2=\bar\sigma(\lambda)=\mu_2-\ii\mu_3\,.
\eeq
Moreover, we consider a basis of eigenvectors of $T$,
also writing the two complex ones in terms of real and imaginary parts
(thus, we do not work directly with complex vectors):
\beq\label{eq:basisT}
  \omega,
  \quad v_2+\ii v_3=\sigma(\omega),
  \quad v_2-\ii v_3=\bar\sigma(\omega),
\eeq
with associated eigenvalues $\lambda$, $\lambda_2$, $\ol\lambda_2$,
respectively.
We understand that, for vectors, the conjugacies $\sigma$, $\bar\sigma$
can be applied componentwisely,
and hence the conjugate vectors above can be obtained just by replacing
$\Omega$ by $\Omega_2$ or $\ol\Omega_2$ in~(\ref{eq:omega}).
In this way, the vectors $v_2$ and $v_3$ do not depend on the specific choice
of a Koch's matrix $T$.
Let $C$ denote the $(3\times3)$-matrix having $\omega$, $v_2$, $v_3$ as columns,
and we consider its condition number
\beq\label{eq:kappa}
  \kappa=\kappa(\omega):=\abs{C}\cdot\abs{C^{-1}},
\eeq
also not depending on the choice of $T$
(we use the matrix norm subordinate to the Euclidean norm for vectors).
Next, we prove that the eigenvalue $\lambda>1$
cannot be arbitrarily close to 1.

\begin{lemma}\label{lm:lambda}
For any Koch's matrix $T\in\SL(3,\Z)$ for $\omega$,
the real eigenvalue $\lambda$ in~(\ref{eq:lambda})
satifies the lower bound $\lambda>\lambda_0$,
with $\lambda_0=\lambda_0(\omega)>1$ defined as the unique real number
satifying $\lambda_0^{\,3}-\lambda_0^{\,2}-\gamma/4\kappa^2=0$,
where $\gamma$ is the constant in the
Diophantine condition~(\ref{eq:DiophCond}),
and $\kappa$ is the condition number~(\ref{eq:kappa}).
\end{lemma}

\proof
From the definitions of $v_2$ and $v_3$, it is clear that
$Tv_2=\mu_2v_2-\mu_3v_3$ and $Tv_3=\mu_3v_2+\mu_2v_3$, and hence
$T=CDC^{-1}$, where we define
$D=\p{\begin{array}{ccc}
         \lambda & 0     &0    \\
         0       & \mu_2 &\mu_3\\
         0       &-\mu_3 &\mu_2\\
       \end{array}}
$.
Since $D\tp D=\diag(\lambda^2,\mu_2^{\,2}+\mu_3^{\,2},\mu_2^{\,2}+\mu_3^{\,2})$,
and using the inequalities
$\sqrt{\mu_2^{\,2}+\mu_3^{\,2}}=\abs{\lambda_2}<1<\lambda$,
one readily sees that $\abs D=\lambda$
and we deduce that $\lambda\le\abs T\le\kappa\lambda$.
Now, we use~(\ref{eq:lambda}), and apply
the Diophantine condition~(\ref{eq:DiophCond}) to the vector
$k=T_{(1)}-(1,0,0)=(T_{11}-1,T_{12},T_{13})$
(it is clear that $k\ne0$, otherwise $T$ has an integer eigenvalue):
\[
  \lambda-1=\scprod k\omega
  \ge\frac{\gamma}{\abs k^2}
  \ge\frac{\gamma}{4\abs T^2}
  \ge\frac{\gamma}{4\kappa^2\lambda^2}\,,
\]
where we used that
$\abs k\le\abs{T_{(1)}}+\abs{(1,0,0)}\le\abs T+1\le2\abs T$.
Finally, a simple study of the function $g(x)=x^3-x^2-\gamma/4\kappa^2$
shows that $\lambda>\lambda_0$.
\qed

Using this lemma, we next show the ``uniqueness'' of the matrix $T$
satisfying Koch's result. More precisely, we can choose
$T=T(\omega)\in\SL(3,\Z)$ whose real eigenvalue $\lambda=\lambda(\omega)>1$
is minimal or, equivalently, the norm $\abs T$ is minimal.
We call this matrix $T$ ``\emph{the principal Koch's matrix for $\omega$}\,''.

\begin{proposition}\label{prop:uniqueT}
There exists a unique matrix $T=T(\omega)\in\SL(3,\Z)$
such that all Koch's matrices for $\omega$ have the form~$\pm T^n$, $n\ge1$.
\end{proposition}

\proof
As we said before, we can restrict ourselves to Koch's matrices of positive
determinant. Assume that $T$ and~$S$ are two Koch's matrices,
with real eigenvalues satisfying $1<\lambda_T\le\lambda_S$.
It is clear that $ST^{-1}$ has $\omega$ as an eigenvector
with eigenvalue $\lambda_S/\lambda_T\ge1$,
and hence $>1$ (it cannot be equal to~1).
This says that $ST^{-1}$ is another Koch's matrix, with
$\lambda_S/\lambda_T>\lambda_0$ by Lemma~\ref{lm:lambda}
(recall that $\lambda_0=\lambda_0(\omega)>1$).
Therefore, the real eigenvalues of the Koch's matrices
for~$\omega$ are all different, and separated at least by a factor $\lambda_0$
(filling in this way a discrete set).
On the other hand, such eigenvalues satisfy the lower bound given
in Lemma~\ref{lm:lambda}. This implies that we can choose a Koch's matrix
$T=T(\omega)$ with minimal eigenvalue $\lambda=\lambda(\omega)>1$.
Then, the matrices $T^n$ (and the opposite ones $-T^n$), $n\ge1$,
are also clearly Koch's matrices.
It remains to show that they are the only ones.
Indeed, if there exists another Koch's matrix $S$,
its real eigenvalue satisfies $\lambda^n<\lambda_S<\lambda^{n+1}$
for some $n\ge1$, and we deduce that  $ST^{-n}$ is a Koch matrix
whose eigenvalue satisfies $1<\lambda_S\lambda^{-n}<\lambda$,
which contradicts our choice of $T$.
\qed

Now, our aim is to describe a simple \emph{procedure} allowing us to determine
the principal Koch's matrix for a given cubic vector $\omega$.
The idea of our method is that any matrix $T$ with integer (or rational)
entries having $\omega$ as an eigenvector is determined
by its first row $T_{(1)}=(T_{11},T_{12},T_{13})$.
The matrices $T$ obtained in this way belong to the general linear group
$\GL(3,\Q)$ but, in general, do not belong to $\SL(3,\Z)$.
However, we can can explore such matrices by giving successive values
to the entries of $T_{(1)}$, until we find a Koch's matrix.
First, in the next lemma we establish the (linear) dependence of $T$
with respect to its first~row.

\begin{lemma}\label{lm:koch}
For any vector $T_{(1)}=(T_{11},T_{12},T_{13})$ with rational
entries, there exists a unique matrix~$T$ with rational entries,
having $\omega$ as an eigenvector, and $T_{(1)}$ as the first row.
This matrix can be written as
\beq\label{eq:T}
  T=A\p{T_{11}\,\Id+T_{12}R+T_{13}(a_0\,\Id+a_1R+a_2R^2)}A^{-1},
\eeq
where we define
\beq\label{eq:R}
  R:=\p{\begin{array}{ccc}
          0   &1   &0\\
          0   &0   &1\\
          r_0 &r_1 &r_2
        \end{array}}
\eeq
(recall the coefficients $r_j$, $a_j$ and the matrix $A$,
introduced in~(\ref{eq:rootOmega}--\ref{eq:tlOmega}) and~(\ref{eq:A})).
\end{lemma}

\proof
We begin by proving the result for the particular case of a frequency vector
$\omega^{(0)}$ as in~(\ref{eq:omega0}).
It is straightforward to check that the matrix $R$
has $\omega^{(0)}$ as an eigenvector with eigenvalue $\Omega$.
The matrix $R^2$, which has $(0,0,1)$ has the first row,
also has the same eigenvector $\omega^{(0)}$ with eigenvalue $\Omega^2$.
Then, it is clear that, for any given vector
$T^{(0)}_{(1)}=\p{T^{(0)}_{11},T^{(0)}_{12},T^{(0)}_{13}}$,
the matrix
\beq\label{eq:T0}
  T^{(0)}=T^{(0)}_{11}\Id+T^{(0)}_{12}R+T^{(0)}_{13}R^2
\eeq
has $T^{(0)}_{(1)}$ as the first row, and $\omega^{(0)}$ as an eigenvector
with eigenvalue
\[
  \lambda=\scprod{T^{(0)}_{(1)}}{\omega^{(0)}}
  =T^{(0)}_{11}+T^{(0)}_{12}\Omega+T^{(0)}_{13}\Omega^2.
\]
To show the uniqueness of such a matrix, notice that its second and third rows
$T^{(0)}_{(2)}$ and $T^{(0)}_{(3)}$ can be determined by the first one
using the equalities $\lambda\Omega=\scprod{T^{(0)}_{(2)}}{\omega^{(0)}}$
and $\lambda\Omega^2=\scprod{T^{(0)}_{(3)}}{\omega^{(0)}}$.
which allow us to determine their entries as (rational) coefficients
in the basis $1$, $\Omega$, $\Omega^2$ of the field~$\Q(\Omega)$.
This shows the result for the particular case of a vector~$\omega^{(0)}$.

Now, we consider the general case of a frequency vector $\omega=A\omega^{(0)}$,
with a matrix $A$ as in~(\ref{eq:A}).
If a matrix $T$ has $\omega$ as an eigenvector
and $T_{(1)}=(T_{11},T_{12},T_{13})$ as the first row,
then it has the form $T=A\,T^{(0)}A^{-1}$,
where $T^{(0)}$ has $\omega^{(0)}$ as an eigenvector,
with the same eigenvalue
\[
  \scprod{T^{(0)}_{(1)}}{\omega^{(0)}}=\lambda
  =\scprod{T_{(1)}}\omega=\scprod{A\tp T_{(1)}}{\omega^{(0)}}
\]
(recall that we consider the rows as column vectors).
Using again that the entries of the vectors can be determined
as coefficients in the basis $1$, $\Omega$, $\Omega^2$,
we deduce that
\[
  T^{(0)}_{(1)}=A\tp T_{(1)}=(T_{11}+a_0T_{13},T_{12}+a_1T_{13},a_2T_{13}).
\]
Applying~(\ref{eq:T0}), we get the whole matrix $T^{(0)}$
and, performing the linear change given by $A$, we get $T$ as in~(\ref{eq:T}).
Its~uniqueness is a direct consequence of the uniqueness of $T^{(0)}$.
\qed

Now, in order to determine the principal Koch's matrix for $\omega$
we can carry out the following simple exploration.
We consider the (integer) entries of the first row $T_{(1)}$
as successive data, say with increasing norm~$\abs{T_{(1)}}$,
until the whole matrix $T$ determined from Lemma~\ref{lm:koch}
belongs to $\SL(3,\Z)$ (i.e.~it has integer entries and determinant~1)
and has an eigenvalue $\lambda>1$ in~(\ref{eq:lambda}).
By Koch's result, we know that such a matrix exists
and will be reached after a finite exploration.
It remains to check whether the matrix $T^*$ obtained in this way
is the principal Koch's matrix for $\omega$ since, in principle,
there could exist another Koch's matrix $T$
with $\abs{T_{(1)}}\ge\abs{T^*_{(1)}}$ but $\abs T<\abs{T^*}$.
If this happens, such a new matrix $T$ would satisfy $\abs{T_{(1)}}<\abs{T^*}$.
Hence, after obtaining a first matrix $T^*$, it is enough to continue
the exploration with increasing norms $\abs{T_{(1)}}$
up to the value $\abs{T^*}$ and, if a new Koch's matrix $T$ is obtained,
check if its norm $\abs T$ is lower than $\abs{T^*}$,
which would imply that the matrix $T$ has to replace $T^*$
as the principal one.

\bremark
In some particular cases, one can provide directly the matrix $A\,R\,A^{-1}$
or its inverse $A\,R^{-1}A^{-1}$ as a Koch matrix. This will happen
if the coefficients $r_j$ and $a_j$ introduced
in~(\ref{eq:rootOmega}--\ref{eq:tlOmega})
are all integer, and $\abs{r_0}=\abs{a_2}=1$.
Since $\det R=r_0$ and $\det A=a_2$, both of the matrices given above
are unimodular (with integer entries and determinant $\pm1$).
Moreover, they have $\omega$ as eigenvector,
with eigenvalue $\Omega$ or $\Omega^{-1}$, respectively.
Notice also that $\Omega$ and $r_0$ have the same sign
(indeed, this comes from the fact that the other two eigenvalues
$\Omega_2$, $\ol\Omega_2$ of $R$ are complex,
and $r_0=\Omega\cdot\Omega_2\cdot\ol\Omega_2$).
We~deduce:
\btm
\item if $\abs\Omega>1$, the matrix $T=r_0A\,R\,A^{-1}$ is a Koch's matrix,
with the eigenvalue $\lambda=r_0\Omega>1$;
\item if $\abs\Omega<1$, the matrix
$T=r_0A\,R^{-1}A^{-1}=-A(r_1\Id+r_2R-R^2)A^{-1}$
is a Koch's matrix, with the eigenvalue $\lambda=r_0\Omega^{-1}>1$.
\etm
However, the Koch's matrix obtained in this way might not be the principal one,
and hence the exploration described above, using the matrices $T$ given
by Lemma~\ref{lm:koch}, would be necessary also in this case.
\eremark

See also in Section~\ref{sect:cubicgolden} the concrete application
of the procedure described above (including the remark)
to the case of the cubic golden vector.
We also recall here that a more systematic algorithm
was developed in \cite{DelshamsGG16} for the case of a
quadratic 2-dimensional vector $\omega=(1,\Omega)$,
providing a ($2\times2$)-matrix $T$, from the (eventually periodic)
continued fraction of the frequency ratio $\Omega$.

Thus, in view of Proposition~\ref{prop:uniqueT},
we will always assume that $T=T(\omega)$ is the principal Koch's matrix.
Since $\det T=1$, it is clear that the modulus the two conjugate eigenvalues is
$\abs{\lambda_2}=\abs{\ol\lambda_2}=\lambda^{-1/2}$.
We now define the following important number,
\beq\label{eq:phi}
  \phi=\phi(\omega)\,:=\,\frac1\pi\,\arg(\lambda_2),
  \qquad \textrm{i.e.} \quad
  \lambda_2,\ol\lambda_2=\frac1{\sqrt\lambda}\,\ee^{\pm\ii\pi\cdot\phi},
\eeq
and we can assume that it is positive: $0<\phi<1$.
Indeed, once the matrix $T(\omega)$ is chosen as the principal one,
the sign of $\phi$ (or equivalently the sign on $\mu_3$ in~(\ref{eq:lambda2}))
is determined by the suitable choice of the sign $s$ for $\sigma_3$
in~(\ref{eq:sigma23b}).

The next lemma has a crucial role in showing that the function $h_1(\eps)$,
appearing in the exponent of the maximal splitting distance
in Theorem~\ref{thm:main},
is quasiperiodic, and not periodic, with respect to $\ln\eps$.
This comes from the fact that the ratio between the two frequencies
of $h_1(\eps)$ is given by $\phi$, as we show in Section~\ref{sect:h1}.

\begin{lemma}\label{lm:irrational}
The number $\phi=\phi(\omega)$ is irrational.
\end{lemma}

\proof
Let us assume that $\phi$ is rational, say $\phi=m/n$ as an irreducible
fraction. Then, the matrix $T^n$ also satisfies Koch's result,
but it has $\lambda^n$ as a simple eigenvalue, and $(-1)^m\lambda^{-n/2}$
as a double real eigenvalue, which contradicts two facts:
the eigenvalues of $T^n$ are all simple, and two of them are complex.
\qed

\subsection{Quasi-resonances of a cubic frequency vector}
\label{sect:resonantseq2}

The matrix $T$ given by Koch's result \cite{Koch99}
provides approximations to the direction of $\omega=(1,\Omega,\wtl\Omega)$.
However, we are not interested in finding approximations to~$\omega$
but, on the contrary, approximations to the quasi-resonances of $\omega$,
which lie close to the  ``resonant plane'' $\langle\omega\rangle^\bot$
(the orthogonal plane to $\omega$).
To be more precise, we say that an integer vector $k\in\Z^3\setminus\pp0$
is a \emph{quasi-resonance} of $\omega$ if
\[
  \abs{\scprod k\omega}<\frac12\,,
\]
and we denote by $\A$ the set of quasi-resonances.

For any given number $x\in\R$, we denote $\rint(x)$ and $\norm x$
the closest integer to $x$ and the distance from $x$ to such closest integer,
respectively.
It is clear that $\norm x=\abs{x-\rint(x)}=\ds\min_{p\in\Z}\abs{x-p}$.
Since the first component of $\omega$ is equal to 1,
for any quasi-resonance $k=(k_1,k_2,k_3)\in\A$
we have $\rint(k_2\Omega+k_3\wtl\Omega)=-k_1$.
In other words, for any $q\in\Z^2\setminus\pp0$ we have a quasi-resonance
\beq\label{eq:k0q}
  k^0(q):=(-p,q)=(-p,q_1,q_2),
  \qquad \textrm{with} \quad p=p^0(q):=\rint(q_1\Omega+q_2\wtl\Omega),
\eeq
whose \emph{small divisor} is
\beq\label{eq:rq}
  r_q:=\scprod{k^0(q)}\omega
  =-p+q_1\Omega+q_2\wtl\Omega
  =\norm{q_1\Omega+q_2\wtl\Omega}.
\eeq
We also say that $k^0(q)$ is an \emph{essential quasi-resonance}
if it is not a multiple of another integer vector,
and we denote by $\A_0$ the set of essential quasi-resonances.

Now, we define the matrix
\beq\label{eq:U}
  U:=(T^{-1})\tp,
\eeq
which satisfies the following simple but important equality:
\beq\label{eq:equalityU}
  \scprod{Uk}\omega=\scprod k{U\tp\omega}
  =\frac1\lambda\scprod k\omega.
\eeq
where $\lambda=\lambda(\omega)$ is the eigenvalue of $T$ with $\lambda>1$.
This says that successive iterations $U^nk$ from a given integer vector~$k$
get closer and closer to the resonant plane~$\langle\omega\rangle^\bot$.

We deduce from~(\ref{eq:equalityU}) that
if $k\in\A$, then also $Uk\in\A$. We say that the
vector $k$ is \emph{primitive} if $k\in\A$ but $U^{-1}k\notin\A$.
It is clear that $k$ is primitive if and only if
the following \emph{fundamental property} is fulfilled:
\beq\label{eq:primitive}
  \frac1{2\lambda}<\abs{\scprod k\omega}<\frac12\,.
\eeq
Writing $k=k^0(q)=(-p,q)$, we denote by $\Pc$ the set of vectors
$q=(q_1,q_2)\in\Z^2\setminus\pp0$ associated to primitive vectors:
\beq\label{eq:calP}
  \Pc:=\{q\in\Z^2: (q_1\ge1\ \textrm{or}\ (q_1=0,q_2\ge1))\ \textrm{and}
  \ k^0(q)\ \textrm{is primitive}\},
\eeq
where the choice of $q_1$ being positive allows us to avoid repetitions,
since it means that $k^0(q)\in\Zc$ (recall the definition~(\ref{eq:calZ})).
We also denote by $\Pc_0$ the set of vectors $q\in\Pc$ such
that $k^0(q)$ is essential.

Now we define, for each $q\in\Pc$,
a \emph{resonant sequence} of integer vectors:
\beq\label{eq:sqn}
  s(q,n) := U^n k^0(q),
  \qquad n\ge0.
\eeq
By construction, the set of such resonant sequences covers the whole set of
quasi-resonances $\A$, providing a classification for them.
As done in \cite{DelshamsG03,DelshamsGG16} for the case of quadratic
frequencies, we are going to establish the properties
of the resonant sequences~(\ref{eq:sqn}) for cubic frequencies
(see Proposition~\ref{prop:cubicfreq} below).

\bremark
A resonant sequence $s(q,n)$ generated by an essential primitive $k^0(q)$
cannot be a multiple of another resonant sequence.
Indeed, in this case we would have $k^0(q)=c\,s(\tl q,n_0)$
with $\abs c>1$ and $n_0\ge0$,
and hence $k^0(q)$ would not be essential.
\eremark

Analogously to the basis of eigenvectors $\omega$, $v_2\pm\ii v_3$ of $T$
introduced in~(\ref{eq:basisT}),
we also consider a basis of eigenvectors of $U$
writing the complex ones in terms of real and imaginary parts:
\beq\label{eq:basisU}
  u_1,
  \quad u_2+\ii u_3=\sigma(u_1),
  \quad u_2-\ii u_3=\bar\sigma(u_1),
\eeq
with eigenvalues
$\lambda^{-1}$, $\lambda_2^{\,-1}$, $\ol\lambda_2^{\,-1}$, respectively.
One readily sees that $\scprod{u_2}\omega=\scprod{u_3}\omega=0$,
i.e.~$u_2$~and~$u_3$ span the resonant plane $\langle\omega\rangle^\bot$.
Other useful equalities are:
$\scprod{u_1}{v_2}=\scprod{u_1}{v_3}=0$,
$\scprod{u_2}{v_2}=-\scprod{u_3}{v_3}$,
$\scprod{u_2}{v_3}=\scprod{u_3}{v_2}$.
We define $Z_1$, $Z_2$ and $\theta$ through the formulas
\beq\label{eq:Z12}
  \frac12(\abs{u_2}^2+\abs{u_3}^2)=Z_1,
  \qquad
  \frac12(\abs{u_2}^2-\abs{u_3}^2)=Z_2\cos\theta,
  \qquad
  \scprod{u_2}{u_3}=Z_2\sin\theta,
\eeq
and the following important number,
\beq\label{eq:delta}
  \delta=\delta(\omega)\,:=\,\frac{Z_2}{Z_1}\,.
\eeq
It is clear, from the definition of $Z_1$ and $Z_2$, that $0\le\delta\le1$.
The following result shows that $\delta$ cannot achieve
the extreme values 0 and 1.
In particular, the fact that $\delta>0$ has a crucial role
(together with the irrationality of $\phi$ shown in Lemma~\ref{lm:irrational})
in showing that the quasiperiodic function $h_1(\eps)$,
appearing in the exponent of the maximal splitting distance
in Theorem~\ref{thm:main}, is not periodic with respect to $\ln\eps$.

\begin{lemma}\label{lm:delta}
The number $\delta=\delta(\omega)$ satisfies $0<\delta<1$.
\end{lemma}

\proof
We first show that $\delta<1$. Indeed, if $\delta=1$ then $Z_1=Z_2$,
which would imply that $\abs{\scprod{u_2}{u_3}}=\abs{u_2}\cdot\abs{u_3}$,
but this is not possible since $u_2$ and $u_3$ are linearly independent.

Now, we are going to see that $\delta>0$. If we have $\delta=0$,
then $Z_2=0$ and, from~(\ref{eq:Z12}),
the expressions $\abs{u_2}^2-\abs{u_3}^2$ and $\scprod{u_2}{u_3}$
would vanish simultaneously.
To show that this is not possible, we are going to see that
they can be written as follows,
\beq\label{eq:u23}
  \abs{u_2}^2-\abs{u_3}^2=c_0+c_1\Omega+c_2\Omega^2,
  \qquad
  \scprod{u_2}{u_3}=(d_0+d_1\Omega+d_2\Omega^2)\,\sigma_3
\eeq
(see~(\ref{eq:sigma23a}) for $\sigma_3$)
and that the coefficients $c_j$, $d_j$ cannot be all zero.

Let us write the coefficients $c_j$, $d_j$ as rational expressions
in the coefficients $r_j$, $a_j$
introduced in~(\ref{eq:rootOmega}--\ref{eq:tlOmega}).
Recall that, in~(\ref{eq:basisU}),
we introduced $u_2\pm\ii u_3$ as complex eigenvectors of the matrix $U$,
conjugates of the real eigenvector~$u_1$.
It is clear from~(\ref{eq:U}) that the eigenvectors of $U$ 
are the same as for $T\tp$.
Since the matrix $T$ can be written as in~(\ref{eq:T})
(with suitable coefficients $T_{1j}$),
it is easy to relate the eigenvectors of $T\tp$ with the ones of $R\tp$,
through the linear change defined by the matrix $B:=(A^{-1})\tp$,
where $A$ is the matrix introduced in~(\ref{eq:A}). Namely, we have
\[
  u_1=B\,u_1^{(0)}, \quad u_2\pm\ii u_3=B\,(u_2^{(0)}\pm\ii u_3^{(0)}),
\]
where
$u_1^{(0)}$, $u_2^{(0)}\pm\ii u_3^{(0)}=\sigma(u_1^{(0)}),\bar\sigma(u_1^{(0)})$
are the eigenvectors of $R\tp$. Using~(\ref{eq:R})
and the cubic equation~(\ref{eq:rootOmega}), it is not hard to obtain the
real eigenvector $u_1^{(0)}$ (with eigenvalue $\Omega$) and, subsequently,
the complex eigenvectors $u_2^{(0)}\pm\ii u_3^{(0)}$ as its conjugates
(with eigenvalues $\sigma_2\pm\ii\sigma_3$, recall~(\ref{eq:sigma23a})).
We get
\bea
  \nonumber
  &&u_1^{(0)}=(r_0,-r_2\Omega+\Omega^2,\Omega),
\\
  \label{eq:u23b}
  &&u_2^{(0)}=(r_0,-r_2\sigma_2+\sigma_2^{\,2}-\sigma_3^{\,2},\sigma_2),
  \quad
  u_3^{(0)}=\sigma_3(0,-r_2+2\sigma_2,1).
\eea
Using such ingredients, together with~(\ref{eq:sigma23b}),
we are able to obtain algebraic expressions for~(\ref{eq:u23})
in the basis $1$, $\Omega$, $\Omega^2$ of the field~$\Q(\Omega)$.
After some tedious computations, we get the following coefficients:
\bean
  &&c_0
  = r_0^{\,2} - \p{\frac{a_0}{a_2}+\frac12}r_0r_2-\frac{2a_1}{a_2}\,r_0
    +r_1^{\,2}-\frac{a_1}{a_2}\,r_1r_2
\\
  &&\phantom{c_0=}
    +\frac{a_0^{\,2}+a_1^{\,2}+1}{a_2^{\,2}}\p{r_1+\frac{r_2^{\,2}}2}\,,
\\
  &&c_1 = \p{\frac{a_0}{a_2}-\frac12}r_0 +\p{\frac{r_2}2+\frac{a_1}{a_2}}r_1\,,
  \qquad
  c_2 = -\frac{r_1}2-\frac{a_0^{\,2}+a_1^{\,2}+1}{2a_2^{\,2}}\,,
\\
  &&d_0=-(c_1+r_2c_2)\,,
  \qquad
  d_1=c_2\,,
  \qquad
  d_2=0\,.
\\
\eean
Assuming $c_j=d_j=0$, $j=0,1,2$, we reach a contradiction.
Indeed, from $c_2=0$ we get $r_1=-(a_0^{\,2}+a_1^{\,2}+1)/a_2^{\,2}$ and,
replacing into the remaining coefficients, we obtain
\bean
  &&c_1=-d_0
    = \p{\frac{a_0}{a_2}-\frac12}r_0
      -\frac{a_0^{\,2}+a_1^{\,2}+1}{a_2^{\,2}}
       \p{\frac{r_2}2+\frac{a_1}{a_2}}\,,
\\
  &&c_0+c_1r_2=r_0^{\,2}-\p{r_2+\frac{2a_1}{a_2}}r_0\,.
\eean
Since $r_0\ne0$ in~(\ref{eq:rootOmega}), from the second equality we get
$r_0=r_2+2a_1/a_2$ and the first equality becomes
\[
  c_1 = -\p{\p{\frac{a_0}{a_2}-1}^2 +\frac{a_1^{\,2}+1}{a_2^{\,2}}}\frac{r_0}2\,,
\]
which contradicts our assumption that $c_1=0$ and,
consequently, we have $\delta>0$.
\qed

\bremark
The previous arguments show that, for the numbers defined
in~(\ref{eq:Z12}), we have $Z_1,Z_2^{\,2}\in\Q(\omega)$.
Indeed, using the rational expressions obtained for the coefficients
$c_j$, $d_j$ (together with the fact that $\sigma_3^{\,2}\in\Q(\omega)$),
we can determine from~(\ref{eq:u23}) the coefficients of $Z_2^{\,2}$
in the basis $1$, $\Omega$, $\Omega^2$.
In an analogous way, we can determine the coefficients of $Z_1$
in the same basis, and we deduce from~(\ref{eq:delta})
that $\delta^2\in\Q(\Omega)$.
Then, it is also possible obtain the coefficients of $\delta^2$
in the basis $1$,~$\Omega$, $\Omega^2$
by carrying out a quotient in the field $\Q(\Omega)$,
though the general expression is very complicated.
See~(\ref{eq:deltagolden}) for the particular case of the
cubic golden frequency vector.
\eremark

For any $q\in\Pc$, we define
\beq\label{eq:gammaq0}
  y_q:=\scprod{k^0(q)}{v_2},
  \qquad
  z_q:=\scprod{k^0(q)}{v_3},
\eeq
and $E_q$, $\psi_q$, $K_q$ and $\gamma^*_q$ through the formulas
\bea
  \label{eq:gammaq1}
  &&\dfrac{\scprod{v_2}{u_2}y_q+\scprod{v_2}{u_3}z_q}
    {\scprod{v_2}{u_2}^2+\scprod{v_2}{u_3}^2}
  =E_q\cos\psi_q\,,
  \quad\!
  \dfrac{\scprod{v_2}{u_3}y_q-\scprod{v_2}{u_2}z_q}
    {\scprod{v_2}{u_2}^2 + \scprod{v_2}{u_3}^2}
  =E_q\sin\psi_q\,,
\\[4pt]
  \label{eq:gammaq2}
  &&K_q:=E_q^{\,2}Z_1\,,
  \qquad
  \gamma^*_q:=\abs{r_q}K_q\,.
\eea

We see in the next proposition that any given resonant sequence $s(q,n)$
defined in~(\ref{eq:sqn}) exhibits
an \emph{``oscillatory limit behavior''} as~$n\to\infty$:
the sizes of the vectors $s(q,n)$ oscillate
around a sequence having geometric growth of rate~$\lambda^{1/2}$,
and the numerators $\gamma_{s(q,n)}$ oscillate around the value $\gamma^*_q$,
which can be considered as the ``mean Diophantine constant''
for the resonant sequence $s(q,n)$.
This proposition extends the results given in
\cite{DelshamsG03,DelshamsGG16} for the quadratic case,
where a (non-oscillatory) limit behavior is also established
for resonant sequences.
In our case of a non-totally real complex vector $\omega$,
the relative amplitude and the frequency of the oscillations
are directly related to the numbers
$\phi=\phi(\omega)$ and $\delta=\delta(\omega)$,
introduced in~(\ref{eq:phi}) and~(\ref{eq:delta}) respectively.
As we see in Section~\ref{sect:asympt_est},
the facts that $\phi$ is irrational and $\delta>0$,
shown by Lemmas~\ref{lm:irrational} and~\ref{lm:delta} respectively,
allow us to show that the function $h_1(\eps)$
associated to the maximal splitting distance in Theorem~\ref{thm:main},
is quasiperiodic but not periodic with respect to $\ln\eps$.

\begin{proposition}\label{prop:cubicfreq}
Let $\omega=(1,\Omega,\wtl\Omega)$ be a cubic frequency vector of complex type.
Consider $\phi$, $\theta$ and $\delta$ as
defined in~(\ref{eq:phi}) and~(\ref{eq:Z12}--\ref{eq:delta}),
and the vector $u_1$ as in~(\ref{eq:basisU}).
For any given $q\in\Pc$,
consider $r_q$, $\psi_q$, $K_q$ and $\gamma^*_q$ as
defined in~(\ref{eq:rq}) and~(\ref{eq:gammaq1}--\ref{eq:gammaq2}).
Then, the resonant sequence $s(q,\cdot)$ defined in~(\ref{eq:sqn})
and its associated numerators $\gamma_{s(q,\cdot)}$ satisfy the approximations
\bea
  \label{eq:cubicfreq_a}
  &&\abs{s(q,n)}^2=K_q\,b_{s(q,n)}\cdot\lambda^n+\Ord(\lambda^{-n/2}),
\\[4pt]
  \label{eq:cubicfreq_b}
  &&\gamma_{s(q,n)}=\gamma^*_q\,b_{s(q,n)}+\Ord(\lambda^{-3n/2}),
\eea
with an oscillating factor defined by
\beq\label{eq:bsqn}
  b_{s(q,n)}:=1+\delta\cos(2\pi\cdot n\phi+2\psi_q-\theta),
\eeq
and hence the numerators $\gamma_{s(q,\cdot)}$ oscillate as $n\to\infty$
between the values
\beq\label{eq:gammapmq}
  \gamma^-_q:=\gamma^*_q\,(1-\delta),
  \qquad
  \gamma^+_q:=\gamma^*_q\,(1+\delta).
\eeq
Moreover, we have the lower bound
\beq\label{eq:lowerbound}
  \gamma^-_q
  \ge\frac{1-\delta}{2\lambda(1+\delta)}(\abs q-Q_0)^2,
  \qquad \textrm{provided} \quad
  \abs q\ge Q_0:=\frac{\abs{u_1}}{2\abs{\scprod{u_1}\omega}}.
\eeq

\end{proposition}

For a proof, see \cite{DelshamsGG14a}.

\bremark
We just outline here the main facts leading
to the dominant behaviors~(\ref{eq:cubicfreq_a}--\ref{eq:cubicfreq_b})
described by this proposition,
and show why this result is valid only in the case
of \emph{complex} conjugates. On one hand, for any given resonant sequence,
the size of the vectors $s(q,n)$ increases like $\lambda^{n/2}$ as $n\to\infty$
(with an oscillatory factor), since
the (coincident) modulus of the greatest eigenvalues of
the iteration matrix $U$ is $\lambda^{1/2}$.
On the other hand, the small divisors
$\abs{\scprod{s(q,n)}\omega}$ decrease like $\lambda^{-n}$
according to the equality~(\ref{eq:equalityU}). Therefore, the numerators
$\gamma_{s(q,n)}=\abs{\scprod{s(q,n)}\omega}\cdot\abs{s(q,n)}^2$
become bounded from above and from below.
This fact does not apply to the \emph{totally real} case,
in which the conjugates of a cubic irrational number
have different modulus.
\eremark

As we can see in~(\ref{eq:cubicfreq_b}), the existence of limit of the
sequences $\gamma_{s(q,n)}$, stated in \cite{DelshamsGG16}
for the quadratic case, is replaced in our complex cubic case by an
oscillatory limit behavior, with a lower limit
$\ds\liminf_{n\to\infty}\gamma_{s(q,n)}=\gamma^-_q$
and an upper limit $\ds\limsup_{n\to\infty}\gamma_{s(q,n)}=\gamma^+_q$\,,
introduced in~(\ref{eq:gammapmq}).
Notice that we could give the exact values of such limits
due to the irrationality of the phase $\phi$
appearing in the oscillating factors~(\ref{eq:bsqn}),
stated in Lemma~\ref{lm:irrational}.

As another relevant fact, we stress that the amplitude
of the limit oscillations is proportional to the number $\delta$
introduced in~(\ref{eq:delta}). Since $\delta>0$ by Lemma~\ref{lm:delta},
we can ensure that such oscillations do occur.

An important consequence of the lower bound~(\ref{eq:lowerbound})
is that the minimal value among the values $\gamma^*_q$
is reached for some concrete $\wh q$.
Indeed, the values $\gamma^*_q$
are not increasing in general with respect to $\abs q$,
but the increasing lower bound~(\ref{eq:lowerbound}) implies that
$\ds\lim_{\abs q\to\infty}\gamma^*_q=\infty$,
and one has to check only a finite number of cases in order to detect
a vector~$\wh q$ providing the minimal value among $\gamma^*_q$, $q\in\Pc$.
We define \emph{the primary resonances} as the integer vectors
belonging to the sequence
\beq\label{eq:s0n}
  s_0(n):=s(\wh q,n),
\eeq
and we denote
\beq\label{eq:gamma}
  \gamma^*:=\min_{q\in\Pc}\gamma^*_q=\gamma^*_{\wh q}>0,
\eeq
which can be considered as the ``minimal mean Diophantine constant''.
The fact that $\gamma^*>0$ implies that
any non-totally real cubic frequency vector $\omega$
satisfies the Diophantine condition~(\ref{eq:diophantine})
(with the minimal exponent~2), and we can compute explicitly
the \emph{``asymptotic Diophantine constant''}~(\ref{eq:gammam}):
\beq\label{eq:asymptdioph}
  \liminf_{\abs k\to\infty}\gamma_k
  =\liminf_{n\to\infty}\gamma_{s_0(n)}=\gamma^*(1-\delta)=\gamma^->0.
\eeq
Dividing by~$\gamma^*$, we also introduce
\emph{normalized numerators} and their associated asymptotic values,
to be used in Section~\ref{sect:asympt_est}:
\beq\label{eq:gammanorm}
  \tl\gamma_k:=\frac{\gamma_k}{\gamma^*}\,,
  \qquad
  \tl\gamma^*_q:=\frac{\gamma^*_q}{\gamma^*}\,,
  \qquad
  \tl\gamma^\pm_q:=\frac{\gamma^\pm_q}{\gamma^*}\,,
\eeq
and in this way we get $\tl\gamma^*_{\wh q}=1$ for the primary resonances.

\bremarks
\item\label{rk:oneprimary}
In principle, for some particular cubic frequency vectors $\omega$,
the minimum in~(\ref{eq:gamma}) could be reached
by two or more vectors $q$ and, consequently, there could exist
two or more sequences of primary resonances.
In such a case, we denote by $\wh q$ only one of such vectors~$q$.
\item\label{rk:primary_ess}
Any primitive vector generating a sequence of primary resonances
is essential: $\wh q\in\Pc_0$. Indeed, if $\wh q$ is not essential,
then we have $k^0(\wh q)=c\,s(\ol q,n_0)$ with $\abs c>1$ and $n_0\ge0$,
and therefore $s(\wh q,n)=c\,s(\ol q,n_0+n)$,
which implies by~(\ref{eq:numerators})
that $\gamma^*_{\wh q}=\abs c^3\gamma^*_{\ol q}$,
and the minimum in~(\ref{eq:gamma})
would not be reached for $\wh q$.
\eremarks

We call \emph{secondary resonances} the vectors
belonging to any of the remaining sequences
$s(q,n)$, $q\in\Pc\setminus\{\wh q\}$.
We also consider the second minimum in~(\ref{eq:gamma}):
\beq\label{eq:gammasec}
  \min_{q\in\Pc\setminus\{\wh q\}}\gamma^*_q=\gamma^*_{\whh q}\,,
\eeq
and we can call ``\emph{main secondary resonances}''
the integer vectors in the sequence $s(\whh q,n)$.
It is clear that its associated normalized numerator
satisfies $\tl\gamma^*_{\whh q}\ge1$.

In order to measure the \emph{``separation''}
between the primary and the secondary resonances,
we define the values
\bea
  \label{eq:Jp0}
  &&J^+_0=J^+_0(\omega):=\p{\tl\gamma^+_{\wh q}}^{1/3}=(1+\delta)^{1/3},
\\
  \label{eq:Bm0}
  &&B^-_0=B^-_0(\omega):=\p{\tl\gamma^-_{\whh q}}^{1/3}
  =\p{\tl\gamma^*_{\whh q}}^{1/3}(1-\delta)^{1/3}
\eea
(we included the exponent $1/3$ for convenience,
see Section~\ref{sect:asympt_est}).
To have a clear distinction between primary and secondary resonances
we need the following ``\emph{weak separation condition}'':
\beq\label{eq:weaksep}
  B^-_0\ >\ J^+_0,
\eeq
which says the interval $[\gamma^-_{\wh q},\gamma^+_{\wh q}]$
has no intersection with any other interval
$[\gamma^-_q,\gamma^+_q]$, $q\ne\wh q$
(as happens for the cubic golden vector, see the next section).

\begin{figure}[!b]
  \centering
  \includegraphics[width=0.6\textwidth,angle=-90]{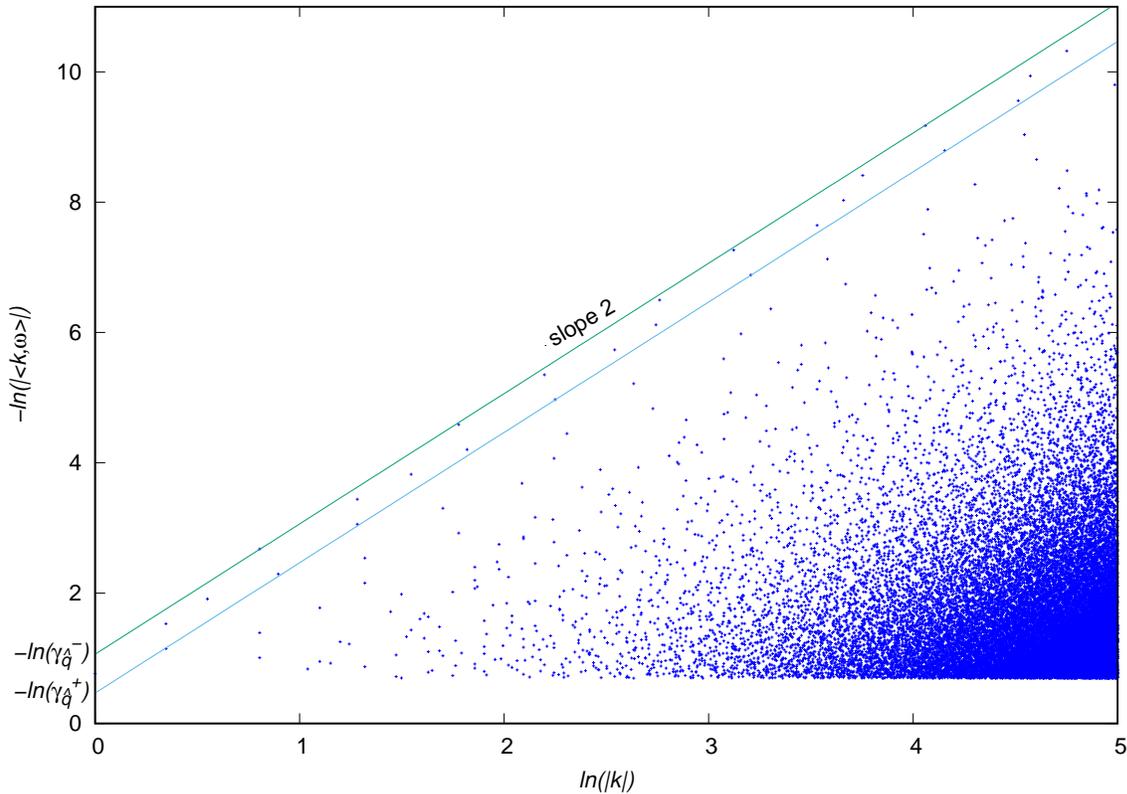}
  \caption{\small\emph{
    Points $(x,y)=(\ln|k|, -\ln|\langle k, \omega\rangle|)$
    for the cubic golden frequency vector;
    the primary resonances correspond to the points
    lying between the two straight lines $y=2x-\ln\gamma^\pm_{\wh q}$.}}
\label{fig:ln_gam}
\end{figure}

Additionally, it is interesting to visualize
the separation between primary and secondary resonances in the following~way.
Taking logarithm of both sides of the Diophantine
condition~(\ref{eq:DiophCond}), we can write it as
\[
  -\ln\abs{\scprod k\omega}\le2\ln\abs k-\ln\gamma.
\]
In Figure~\ref{fig:ln_gam} (which corresponds to the cubic golden vector),
where we draw all the points with coordinates
$(x,y)=(\ln\abs k,-\ln\abs{\scprod k\omega})$
(up to a large value of $\abs k$),
we can see a sequence of points lying between the two straight lines
$y=2x-\ln\gamma^\pm_{\wh q}$.
Those points correspond to integer vectors belonging
to the sequence of primary resonances: $k=s_0(n)$, $n\ge0$,
and the remaining points correspond to secondary resonances.

\subsection{The cubic golden frequency vector}
\label{sect:cubicgolden}

In this section, we provide particular data
for the concrete case of the cubic golden frequency vector.
We point out that a similar approach could be carried out for other
cubic vectors (see \cite{Chandre02} for some famous examples).

We introduce $\Omega$ as the real number satisfying $\Omega^3=1-\Omega$,
which has been called the \emph{cubic golden number}
(see for instance \cite{HardcastleK00}).
Then, we consider the frequency vector
\beq\label{eq:cubicgolden}
  \omega=(1,\Omega,\Omega^2)\approx(\,1\,,\,0.682328\,,\,0.465571\,).
\eeq
In other words, the coefficients introduced
in~(\ref{eq:rootOmega}--\ref{eq:tlOmega})
are $r_0=1$, $r_1=-1$, $r_2=0$, $a_0=a_1=0$, $a_2=1$,
and hence the matrices defined in~(\ref{eq:R}) and~(\ref{eq:A})
are $R=\p{\begin{array}{ccc}
         0 & 1 &0\\
         0 & 0 &1\\
         1 &-1 &0\\
       \end{array}}
$
and $A=\Id$.

In fact, we can provide exact expressions for $\Omega$ using
some of the standard formulas for the solutions
of the general cubic equation (see for instance \cite{Weisstein03}).
We have
\[
  \Omega=S_++S_-=S_\pm-\frac1{3S_\pm}\,,
  \qquad \textrm{with} \quad
  S_\pm=\sqrt[3]{\frac12\p{1\pm\sqrt{\frac{31}{27}}}}\,,
\]
or also
\[
  \Omega=\frac2{\sqrt3}\sinh\p{\frac13\arsinh\frac{3\sqrt3}2}.
\]

It is easy, from the results of Section~\ref{sect:resonantseq1},
to obtain the principal Koch's matrix
for the frequency vector~(\ref{eq:cubicgolden}).
By Lemma~\ref{lm:koch}, any Koch's matrix is determined
from its first row $T_{(1)}=(T_{11},T_{12},T_{13})$,
by the formula $T=T_{11}\,\Id+T_{12}R+T_{13}R^2$.
On the other hand, the remark after Lemma~\ref{lm:koch}
ensures that $T^*=R^{-1}=\Id+R^2$ is a Koch's matrix but, in principle,
it might not be the principal one.
To check whether another Koch's matrix can be the principal one,
we carry out the exploration described after Lemma~\ref{lm:koch}
in the following way.
We use that the matrix $T^*$ given above
has norm $\abs{T^*}=(\sqrt5+1)/2\approx1.618034$,
and its first row $T^*_{(1)}=(1,0,1)$ has norm
$\abs{T^*_{(1)}}=\sqrt2\approx1.414214$.
Then, by exploring the matrices $T$ given by a few possible first rows $T_{(1)}$
(with norms between $\sqrt2$ and $(\sqrt5+1)/2$),
we ensure that the Koch's matrix $T^*$ given above is the principal one.
We rename it as $T$.

In this way, the principal Koch's matrix
for the cubic golden frequency vector~(\ref{eq:cubicgolden}),
and the subsequent matrix introduced in~(\ref{eq:U}), are
\[
  T=R^{-1}=\p{\begin{array}{ccc}
         1 &0 &1\\
         1 &0 &0\\
         0 &1 &0\\
       \end{array}},
  \qquad
  U=R\tp
  =\p{\begin{array}{ccc}
        0 &0 &1\\
        1 &0 &-1\\
        0 &1 &0\\
      \end{array}},
\]
with the eigenvalue
\beq\label{eq:lambdagolden}
  \lambda=1+\Omega^2=\frac1\Omega\approx1.465571\,,
\eeq
which satisfies $\lambda^3=1+\lambda^2$.

Let us compute several relevant parameters, defined
in Section~\ref{sect:resonantseq1}.
Writing the conjugates of $\Omega$
as $\Omega_2,\ol\Omega_2=\sigma_2\pm\ii\sigma_3$,
by~(\ref{eq:sigma23b}) we have
\[
  \sigma_2=-\frac\Omega2\,,
  \qquad
  \sigma_3=-\frac{\ds\sqrt{4+3\Omega^2}}2\,,
\]
where the sign $s=-1$ chosen for $\sigma_3$ in~(\ref{eq:sigma23b})
ensures that $\lambda_2=1/\Omega_2=1/(\sigma_2+\ii\sigma_3)$
has positive imaginary part,
and hence the the number defined in~(\ref{eq:phi}) is
\[
  \phi=1+\frac1\pi\,\arctan\frac{-\sigma_3}{\sigma_2}
  \approx0.590935\,,
\]
and it is irrational by Lemma~\ref{lm:irrational}.
As stated in Theorem~\ref{thm:main}, the number $\phi$
is the frequency ratio of the function $h_1(\eps)$ as a quasiperiodic funtion
(with respect to $\ln\eps$).
It is interesting to consider its (infinite) \emph{continued fraction}
and its associated convergents,
whose denominators provide ``approximate periods'' for $h_1(\eps)=F_1(\zeta)$
(in the logarithmic variable $\zeta\sim\ln(1/\eps)$, see~(\ref{eq:logchange})):
\[
  \phi=[\,0;\,1,1,2,4,78,\ldots\,]
  \approx
  \frac11\,,\frac12\,,\frac35\,,\frac{13}{22}\,,\frac{1017}{1721}\,,\ldots
\]
In particular, the convergent $13/22$ is close enough to $\phi$,
and explains the fact that $F_1(\zeta)$
appears to be 22-periodic in Figure~\ref{fig:h1}.
On the other hand, the number $\delta$ introduced in~(\ref{eq:delta})
can be obtained by carrying out, for this particular case,
the computations described in the remark after Lemma~\ref{lm:delta},
and we get
\beq\label{eq:deltagolden}
  \delta=\sqrt{-1+5\Omega-5\Omega^2}\approx0.289453\,.
\eeq

In the table below, we write down several numerical data
appearing in Proposition~\ref{prop:cubicfreq},
for the resonant sequences $s(q,n)$ induced by the primitives~$k^0(q)$
(see~(\ref{eq:k0q}) and~(\ref{eq:sqn})):
the numbers $\gamma^*_q$, the bounds $\gamma^-_q$ and~$\gamma^+_q$,
and the normalized values $\tl\gamma^*_q$
(defined in~(\ref{eq:gammaq0}--\ref{eq:gammaq2}), (\ref{eq:gammapmq})
and~(\ref{eq:gammanorm}), respectively;
we also use the expressions~(\ref{eq:basisT}) and~(\ref{eq:u23b})
for the vectors $v_j$ and $u_j$).
We restrict such data to the primitives~$k^0(q)$ with $\abs q<3$,
and we provide a lower bound for all other primitives
(see~(\ref{eq:lowerbound})).

\begin{center}
\begin{tabular}{|c||c|c|c||c|}
  \hline
  $k^0(q)=(-p,q)$&$\gamma^-_q$ &$\gamma^*_q$&$\gamma^+_q$&$\tl\gamma^*_q$\\[2pt]
  \hline\hline
  $(0,0,1)$      &0.345858     &0.486749    &0.627640    &1\\
  $(-1,2,0)$     &1.037575     &1.460248    &1.882920    &3\\
  $(-2,1,2)$     &3.112725     &4.380743    &5.648761    &9\\
  $(0,2,-2)$     &2.766867     &3.893994    &5.021121    &8\\[2pt]
  \hline
  $\abs q\ge3$   &$\ge1.274218$&\multicolumn3{c|}{}\\
  \hline
\end{tabular}
\end{center}

As we see from this table, the smallest value of $\gamma^*_q$ corresponds
to $\wh q=(0,1)$, i.e.~to the primitive vector $k^0(\wh q)=(0,0,1)$,
which generates the sequence of primary resonances.
The minimum of the values~$\gamma^*_q$ is the
``minimal mean Diophantine constant'' introduced in~(\ref{eq:gamma}):
\[
  \gamma^*=\gamma^*_{\wh q}=\dfrac2{31}(5+\Omega+4\Omega^2)\approx0.486749
\]
(the algebraic expression in the basis $1$, $\Omega$, $\Omega^2$ has also
been obtained from the definition~(\ref{eq:gammaq0}--\ref{eq:gammaq2}),
working in the field~$\Q(\Omega)$).
On the other hand, we get for the
``asymptotic Diophantine constant''~(\ref{eq:asymptdioph})
the value $\gamma^-\approx0.345858$\,.
Other numerical values appearing in Proposition~\ref{prop:cubicfreq}
are $\theta\approx-1.054837$ and $\psi_{\wh q}\approx-2.007416$
(the latter one for the primary resonances),
defined in~(\ref{eq:Z12}) and~(\ref{eq:gammaq1}) respectively.

Finally, in~(\ref{eq:Jp0}--\ref{eq:Bm0}) we get
\beq\label{eq:JB0golden}
  J^+_0=(1+\delta)^{1/3}\approx1.088433\,,
  \qquad
  B^-_0=3^{1/3}(1-\delta)^{1/3}\approx1.286979\,,
\eeq
and hence the weak separation condition~(\ref{eq:weaksep}) is fulfilled.

\section{Searching for the asymptotic estimate}
\label{sect:asympt_est}

In order to provide an asymptotic estimate for the splitting,
given in our main result (Theorem~\ref{thm:main})
in terms of the splitting function $\M(\theta)$,
we first need to carry out a careful study of the first order
approximation~(\ref{eq:melniapproxM})
provided by the Poincar\'e--Melnikov method.
Although this approximation is given by the (vector) Melnikov function
$M(\theta)$, $\theta\in\T^3$, it is more convenient to work with
the (scalar) splitting potential $L(\theta)$,
whose gradient is the Melnikov function: $\nabla L(\theta)=M(\theta)$.

In this section, we provide the \emph{constructive} part of the proof,
which amounts to find, for every sufficiently small $\eps$,
the dominant harmonic of the Fourier expansion
of the Melnikov potential $L(\theta)$,
with an asymptotic estimate for its size
of the type $\exp\{-h_1(\eps)/\eps^{1/6}\}$,
with an oscillating (positive) function $h_1(\eps)$ in the exponent.
This function can be explicitly defined from the arithmetic properties
of our cubic frequency vector $\omega$ and, as a direct consequence,
we see that it is quasiperiodic (and continuous) with respect to $\ln\eps$,
and hence bounded (and we provide concrete lower and upper bounds for it).
We can also study, from such arithmetic properties,
whether the dominant harmonic is always given by a primary resonance
(providing a sufficient condition for this, which is satisfied in the case
of the cubic golden frequency vector) or, otherwise,
secondary resonances can be dominant for some intervals of $\eps$.

The final step, considered in Section~\ref{sect:technical},
requires to ensure that the whole Melnikov function $M(\theta)$
is dominated by its dominant harmonic,
by obtaining a bound for the sum of all the remaining harmonics
of its Fourier expansion.
Furthermore, to ensure that the Poincar\'e--Melnikov
method~(\ref{eq:melniapproxM})
predicts correctly the size of the splitting in the singular case $\mu=\eps^r$,
one has to extend the results to the splitting function $\M(\theta)$
by showing that the asymptotic estimate of the dominant harmonic
is large enough to overcome the harmonics of the error term
in~(\ref{eq:melniapproxM}).
This step is just outlined in Section~\ref{sect:technical},
since it is analogous to the one already done in \cite{DelshamsG04}
for the case of the quadratic golden number
(using the upper bounds for the error term provided in \cite{DelshamsGS04}).

\subsection{Estimates of the harmonics of the splitting
  potential}
\label{sect:gn}

We plug our functions $f$ and $h$, defined in~(\ref{eq:hf}),
into the integral~(\ref{eq:L}) and get the Fourier expansion
of the Melnikov potential,
where the coefficients can be obtained using residues
(see for instance \cite[\S3.3]{DelshamsG00}):
\beq\label{eq:coefL}
  L(\theta) = \sum_{k\in\Zc\setminus\pp0}
  L_k\,\cos(\langle k, \theta\rangle -\sigma_k),
  \qquad
  L_k = \frac{2\pi |\langle k, \omega_\eps\rangle|
  \,\ee^{-\rho\abs k}}{\sinh |\frac{\pi}{2}
  \langle k, \omega_\eps\rangle|}\,,
\eeq
and the phases $\sigma_k$ are the same as in~(\ref{eq:hf}).
Recalling that the fast frequencies $\omega_\eps$ are given
in~(\ref{eq:omega_eps}) and taking into account the
definition of the numerators $\gamma_k$ in~(\ref{eq:numerators}),
we can present each coefficient $L_k=L_k(\eps)$, $k\in\Zc\setminus\pp0$
(recall that we introduced the set $\Zc\subset\Z^3$ in~(\ref{eq:calZ}),
to avoid repetitions in Fourier expansions),
in the form
\begin{equation}
\label{eq:alphabeta}
  L_k = \alpha_k\,\ee^{- \beta_k},
  \quad
  \beta_k(\eps)=\rho\abs k+\frac{\pi\gamma_k}{2\abs k^2\sqrt{\eps}}\,,
  \qquad
  \alpha_k(\eps)\approx\frac{4\pi\gamma_k}{\abs k^2\sqrt{\eps}}\,,
\end{equation}
where an exponentially small term has been neglected in the denominator
of $\alpha_k$.
The most relevant term in this expression is $\beta_k$, which gives
the exponential smallness in $\eps$ of each coefficient,
and we will show that $\alpha_k$ provides a~polynomial factor.
For any given $\eps$, the smallest exponents $\beta_k(\eps)$
provide the largest (exponentially small) coefficients $L_k(\eps)$
and hence the dominant harmonics. Our aim is to study the dependence
on $\eps$ of the size of the most dominant harmonic.

To start, we provide a more convenient expression
for the exponents $\beta_k(\eps)$,
which shows that the smallest ones are~$\Ord(\eps^{-1/6})$.
Indeed, we deduce from~(\ref{eq:alphabeta}) that we can write
\beq\label{eq:beta_gk}
  \beta_k(\eps) = \frac{C_0}{\eps^{1/6}}\,g_k (\eps),
  \qquad
  C_0:=\frac32(\pi\rho^2\gamma^*)^{1/3},
\eeq
where for any given $k$ we introduce the function
\[
  g_k (\eps):=\frac{\tl\gamma_k^{1/3}}3
  \pq{2\p{\frac\eps{\eps_k}}^{1/6}+\p{\frac{\eps_k}\eps}^{1/3}},
  \qquad
  \eps_k:=\frac{D_0\tl\gamma_k^{\,2}}{\abs k^6}\,,
  \qquad
  D_0:=\p{\frac{\pi\gamma^*}\rho}^2.
\]
It is straightforward to check that each function $g_k(\eps)$
attains its minimum at $\eps=\eps_k$,
with the (positive) minimum value $g_k(\eps_k)=\tl\gamma_k^{\,1/3}$.
Recall that the constant $\gamma^*=\gamma^*_{\wh q}$
and the normalized numerators $\tl\gamma_k=\gamma_k/\gamma^*$
were introduced in~(\ref{eq:gamma}) and~(\ref{eq:gammanorm}), respectively.

Since we are interested in obtaining \emph{asymptotic estimates}
for the splitting distance, rather than lower bounds,
we need to determine for any given $\eps$ the most dominant harmonic,
which is given by the smallest value $g_k(\eps)$,
reached for some integer vector $k=S_1(\eps)$ to be determined.
In fact, as in \cite{DelshamsGG16} we may replace, for $\eps$ small,
the functions $g_k(\eps)$ by approximations $g^*_k(\eps)$,
obtained by neglecting the asymptotic terms
going to 0 in Proposition~\ref{prop:cubicfreq}.
More precisely, for $k=s(q,n)$ belonging to a concrete resonant sequence,
we use the approximations~(\ref{eq:cubicfreq_a}--\ref{eq:cubicfreq_b})
for $\abs{s(q,n)}$ and $\gamma_{s(q,n)}$ as $n\to\infty$,
given in Proposition~\ref{prop:cubicfreq},
and we obtain the following approximations:
\bea
  \label{eq:gqn1}
  &&g_{s(q,n)}(\eps) \approx g^*_{s(q,n)}(\eps)
  :=\frac{(\tl\gamma^*_qb_{s(q,n)})^{1/3}}3
  \pq{2\p{\frac\eps{\eps^*_{s(q,n)}}}^{1/6}
    +\p{\frac{\eps^*_{s(q,n)}}\eps}^{1/3}},
\\
  \label{eq:gqn2}
  &&\eps_{s(q,n)} \approx \eps^*_{s(q,n)}
  :=\frac{D_0(\tl\gamma^*_q)^2}{K_q^{\,3}\,b_{s(q,n)}\cdot\lambda^{3n}}\,,
\eea
with the oscillating factors $b_{s(q,n)}$ introduced in~(\ref{eq:bsqn}).
Notice that each function $g^*_{s(q,n)}(\eps)$ has its minimum
at~$\eps^*_{s(q,n)}$, whose dependence on $n$ is not strictly geometric
(decreasing with ratio $\lambda^3$),
but ``perturbed'' by the oscillating factor~$b_{s(q,n)}$.
Analogously, the minimum values
$g^*_{s(q,n)}(\eps^*_{s(q,n)})=\tl\gamma^*_qb_{s(q,n)}$
are not constant but oscillating.
The size of such ``perturbations'' is given by the value $\delta$
introduced in~(\ref{eq:delta}).

\bremark
The most dominant harmonic cannot be found
in a non-essential resonant sequence.
Indeed, if $s(q,n)=c\,s(\ol q,n_0+n)$ with $\abs c>1$ and $n_0\ge0$, then
$g^*_{s(q,n)}(\eps)=\abs c\,g^*_{s(\ol q,n_0+n)}(\eps)$
(see also remark~\ref{rk:primary_ess}
at the end of Section~\ref{sect:resonantseq2}).
\eremark

The sequence of primary resonances $s_0(n)=s(\wh q,n)$,
defined in~(\ref{eq:s0n}), plays an important role
since it gives the smallest minimum values among the functions $g^*_k(\eps)$,
and hence they will provide the most dominant harmonics,
at least for $\eps$ close to such minima.
With this fact in mind, and recalling that $\tl\gamma^*_{\wh q}=1$, we introduce
\bea
  \label{eq:gn1}
  &&\bg_n(\eps)
  :=g^*_{s_0(n)}(\eps)
  =\frac{\bb_n^{\,1/3}}3\pq{2\p{\frac\eps{\beps_n}}^{1/6}
   +\p{\frac{\beps_n}\eps}^{1/3}},
\\
  \label{eq:gn2}
  &&\beps_n:=\eps^*_{s_0(n)}=\frac{D_0}{K_{\wh q}^{\,3}\bb_n\cdot\lambda^{3n}}\,,
\\
  \label{eq:gn3}
  &&\bb_n:=b_{s_0(n)}=1+\delta\cos(2\pi\cdot n\phi+2\psi_{\wh q}-\theta).
\eea

In order to determine the most dominant harmonic for any given $\eps$,
we have to study the relative position of the functions~$g^*_{s(q,n)}(\eps)$
and the intersections between their graphs.
Due to the (essentially) geometric behavior
of the minima $\eps^*_{s(q,n)}$ as $n\to\infty$,
it is convenient to replace $\eps$ by a logarithmic variable:
\beq\label{eq:logchange}
  \zeta=\Lg\frac{D_0}{K_{\wh q}^{\,3}}-\Lg\eps,
  \qquad \textrm{i.e.} \quad \eps=\frac{D_0}{K_{\wh q}^{\,3}\,\lambda^{3\zeta}}
\eeq
(notice that $\zeta\to\infty$ as $\eps\to0^+$),
where we introduce the notation
\[
  \Lg x:=\log_{(\lambda^3)}x=\frac{\ln x}{3\ln\lambda}\,.
\]
We define for any given $Z\in\R$ and $Y>0$
the following \emph{``hyperbolic cosine-like'' function}:
\beq\label{eq:defC}
  \Cc(\zeta\,;\,Z,Y):=Y^{1/3}\,\Cc_0(\zeta-Z),
  \qquad
  \Cc_0(\zeta):=\frac13(2\lambda^{-\zeta/2}+\lambda^\zeta).
\eeq
Any function $\Cc(\zeta\,;\,Z,Y)$ has its minimum at $\zeta=Z$
with $\Cc(Z\,;\,Z,Y)=Y^{1/3}$ as the minimum value, and is
a \emph{convex} function.
In fact, the point~$(Z,Y^{1/3})$ of its graph determines the function, and
the graph becomes divided at this point
into a ``decreasing branch'' ($\zeta<Z$)
and an ``increasing branch'' ($\zeta>Z$).

Translating definitions~(\ref{eq:gqn1}--\ref{eq:gn2})
of $g^*_{s(q,n)}(\eps)$, $\eps^*_{s(q,n)}$, $\bg_n(\eps)$, $\beps_n$
into the new variable, we get:
\bea
  \label{eq:fqn1}
  &&f^*_{s(q,n)}(\zeta)
  :=\Cc(\zeta\,;\,\zeta^*_{s(q,n)},\tl\gamma^*_qb_{s(q,n)}),
\\
  \label{eq:fqn2}
  &&\zeta^*_{s(q,n)}
  :=n+3\Lg\frac{K_q}{K_{\wh q}}-2\Lg\tl\gamma^*_q+\Lg b_{s(q,n)}\,,
\\
  \label{eq:fqn3}
  &&\bbf_n(\zeta):=\Cc(\zeta\,;\,\bzeta_n,\bb_n),
  \qquad\qquad
  \bzeta_n:=n+\Lg\bb_n\,,
\eea
Notice that, if the oscillating terms $b_{s(q,n)}$ are not taken into account
(i.e.~if we assume $\delta=0$ in~(\ref{eq:delta})),
the graph of a function $f^*_{s(q,n+1)}$ is a translation of $f^*_{s(q,n)}$
to distance~1, which would be the situation for the case of
quadratic frequencies considered in \cite{DelshamsGG16}.
What we actually have for cubic frequencies
is an $\Ord(\delta)$-perturbation of this situation,
due to the terms~$b_{s(q,n)}$ defined in~(\ref{eq:bsqn}).

\bremark
In fact, if analogous computations are carried out for the quadratic case,
the function $\Cc_0(\zeta)$ introduced in~(\ref{eq:defC}) should be
replaced by an expression of the type
$(\lambda^{-\zeta}+\lambda^\zeta)/2=\cosh(\zeta\ln\lambda)$
(with a somewhat different definition of the variable $\zeta$).
An expression of this type in asymptotic estimates for the splitting
appeared for the first time in \cite{DelshamsGJS97}
(see also \cite{DelshamsG04}).
We point out that our ``hyperbolic cosine-like'' function $\Cc_0(\zeta)$
is no longer an even function of $\zeta$ in the cubic case considered here,
according to the definition~(\ref{eq:defC}).
In other words, the symmetry of the ``true'' hyperbolic cosine function
$\cosh(\zeta\ln\lambda)$ between the decreasing and increasing branches,
that takes place in the quadratic case, is not preserved in the cubic case.
\eremark

In order to study the dependence of the most dominant harmonics
on $\eps$, now replaced by
the logarithmic variable~$\zeta$ introduced in~(\ref{eq:logchange}),
it is useful to consider the intersections
between the graphs of functions~(\ref{eq:fqn1}),
since this gives the values of $\zeta$ at which
a change in the dominance may take place.
The next two lemmas show that, if we consider the graphs associated
to the functions $f^*_k(\zeta)$ and $f^*_{\ol k}(\zeta)$
associated to different quasi-resonances $k$, $\ol k$,
only two situations are possible:
they do not intersect (which says that one of them always dominates
the other one), or they intersect transversely at a unique point
(and in this case a unique change in the dominance takes place
among such two quasi-resonances).
Namely, in Lemma~\ref{lm:fkk} we show that $f^*_k$ and $f^*_{\ol k}$
cannot be the same function,
and in Lemma~\ref{lm:ZW}
(formulated, by convenience, in terms of the functions
introduced in~(\ref{eq:defC}))
we provide the condition for the existence
of intersection between their graphs, as well as an explicit formula
for this intersection,
and some additional bounds to be used later.

\begin{lemma}\label{lm:fkk}
For any given $k,\ol k\in\A\cap\Zc$ with $k\ne\ol k$,
the functions $f^*_k(\zeta)$ and $f^*_{\ol k}(\zeta)$ do not coincide.
\end{lemma}

\proof
Recalling the definition~(\ref{eq:sqn}),
let us write $k=s(q,n)$ and  $\ol k=s(\ol q,\ol n)$.
If $f^*_k=f^*_{\ol k}$, then we have $g^*_k=g^*_{\ol k}$ and,
by definition~(\ref{eq:gqn1}), we get
$\tl\gamma^*_qb_k=\tl\gamma^*_{\ol q}b_{\ol k}$ and $\eps^*_k=\eps^*_{\ol k}$\,.
By~(\ref{eq:gammaq2}), such two equalities can be rewritten as
$\abs{r_q}K_qb_k=\abs{r_{\ol q}}K_{\ol q}b_{\ol k}$ and
$K_qb_k\lambda^n=K_{\ol q}b_{\ol k}\lambda^{\ol n}$, respectively.
We deduce that the small divisors~(\ref{eq:rq}) satisfy
$\abs{r_{\ol q}/r_q}=\lambda^{\ol n-n}$ but,
from the fundamental property~(\ref{eq:primitive}),
we have $\abs{r_q},\abs{r_{\ol q}}\in(1/2\lambda,1/2)$.
This says that $n=\ol n$ and hence $\abs{r_q}=\abs{r_{\ol q}}$,
but from definition~(\ref{eq:rq}) and the fact that $\omega$
is a nonresonant vector we deduce that $q=\pm\ol q$,
which contradicts the assumption~$k\ne\ol k$
(recall that $k,\ol k\in\Zc$).
\qed

\begin{lemma}\label{lm:ZW}
Let $Z_1,Z_2\in\R$ and $Y_1,Y_2>0$ with $(Z_1,Y_1)\ne(Z_2,Y_2)$, and define
\[
  Z=Z_2-Z_1,
  \qquad
  W=\p{\frac{Y_2}{Y_1}}^{1/3}.
\]
Then, we have:
\btm
\item[(a)] The graphs of the functions
$\Cc(\zeta\,;\,Z_1,Y_1)$ and $\Cc(\zeta\,;\,Z_2,Y_2)$
intersect if and only if
$\lambda^Z<\min(W,W^{-2})$ or $\lambda^Z>\max(W,W^{-2})$.
If~so, the intersection is unique and transverse,
and takes place at the point given by
\beq\label{eq:ZW}
  \zeta^*=Z_1+2\Lg\frac{2\lambda^Z(W\lambda^{Z/2}-1)}{\lambda^Z-W}\,.
\eeq
\item[(b)] The following upper/lower bound holds:
\[
  \begin{array}{ll}
    \zeta^*<Z_1+2\Lg\dfrac{2\lambda^Z}{W-\lambda^Z}
    &\textrm{if} \quad \lambda^Z<\min(W,W^{-2}),
  \\[6pt]
    \zeta^*>Z_1+2\Lg2(W\lambda^{Z/2}-1)
    &\textrm{if} \quad \lambda^Z>\max(W,W^{-2}).
  \end{array}
\]
\etm
\end{lemma}

\proof
Introducing the variable $\xi=\zeta-Z_1$, we see from
definition~(\ref{eq:defC}) that the intersection between
the graphs of $\Cc(\zeta\,;\,Z_1,Y_1)$ and $\Cc(\zeta\,;\,Z_2,Y_2)$
corresponds to the solution of the equation $\Cc_0(\xi)=W\,\Cc_0(\xi-Z)$,
where we have $(Z,W)\ne(0,1)$.
After some straightforward computations, we see that this solution $\xi=\xi^*$
is given by
$\lambda^{3\xi^*/2}=\dfrac{2\lambda^Z(W\lambda^{Z/2}-1)}{\lambda^Z-W}$\,,
which leads directly to the formula~(\ref{eq:ZW}) for $\zeta^*=Z_1+\xi^*$.
Notice that the intersection does not take place if
$\lambda^Z$ belongs to the interval of endpoints $W$ and $W^{-2}$
(indeed, in this case the numerator and denominator
in the expression~(\ref{eq:ZW}) would have different sign).

To complete the proof of~(a), we have to show the transversality
of the intersection. This amounts to see
that the solution obtained above does not satisfy
the equation $\Cc_0'(\xi^*)=W\,\Cc_0'(\xi^*-Z)$.
Indeed, solving this new equation we get
$\lambda^{3\xi^*/2}=\dfrac{\lambda^Z(W\lambda^{Z/2}-1)}{W-\lambda^Z}$\,,
which is possible only if
$\lambda^Z$ does belong to the interval of endpoints $W$ and $W^{-2}$
(the~case excluded above).

The proof of the bound~(b) for $\zeta^*$, in the two cases considered,
is straightforward from the formula~(\ref{eq:ZW}).
\qed

\subsection{Estimate of the most dominant harmonic}
\label{sect:h1}

We introduce the positive function $h_1(\eps)$
appearing in the exponent in Theorem~\ref{thm:main}
as the minimum, for any given $\eps$,
of the values $g^*_k(\eps)$ among the quasi-resonances,
and we denote $S_1=S_1(\eps)$ the integer vector $k$
at which such minimum is reached:
\beq\label{eq:h1}
  h_1(\eps):=\min_{k\in\A}g^*_k(\eps)=g^*_{S_1}(\eps).
\eeq
In fact, by the remark after definitions~(\ref{eq:gqn1}--\ref{eq:gqn2})
the integer vector providing the minimum is always
an essential quasi-resonance: $S_1(\eps)\in\A_0$.

Our aim is to study some of the properties of $h_1(\eps)$,
putting emphasis on the dependence of such functions
on the arithmetic properties of the cubic frequency vector $\omega$,
studied in Section~\ref{sect:cubicfreq}.
Namely, we prove that the function $h_1(\eps)$ satisfies
the following properties:
\btm
\item It is \emph{piecewise-smooth} and \emph{piecewise-convex}
(and continuous), with corners (i.e.~jump discontinuities of the derivative)
associated to changes in the dominant harmonic
(i.e.~discontinuities of the ``piecewise-constant'' function~$S_1(\eps)$).
\item It is \emph{bounded}, providing (positive) lower and upper bounds
for it.
\item It is \emph{quasiperiodic} (and \emph{not periodic})
with respect to $\ln\eps$, with two frequencies whose ratio
is the irrational number $\phi$ defined in~(\ref{eq:phi}).
\etm

As in Section~\ref{sect:gn}, we can translate the function $h_1(\eps)$
into the logarithmic variable $\zeta$ introduced in~(\ref{eq:logchange}):
\[
  F_1(\zeta):=\min_{k\in\A}f^*_k(\zeta)=f^*_{R_1}(\zeta),
\]
with $R_1=R_1(\zeta)=S_1(\eps)$.
We also define an analogous but somewhat simpler function,
taking into account \emph{only the primary} resonances $s_0(n)$
introduced in~(\ref{eq:s0n}) and involved
in~(\ref{eq:gn3}) and~(\ref{eq:fqn3}):
\beq\label{eq:bF1}
  \bF_1(\zeta):=\min_{n\ge0}\bbf_n(\zeta)=\bbf_{N_1}(\zeta),
\eeq
with $N_1=N_1(\zeta)$. In other words, the most dominant harmonic
among the primary resonances corresponds to $\bR_1=\bR_1(\zeta)=s_0(N_1)$.

Clearly, for any $\zeta$ we have
\beq\label{eq:FF1}
  F_1(\zeta)\le\bF_1(\zeta).
\eeq
In order to provide an accurate description of the splitting,
it is useful to study whether the equality between
the above functions can be established for any value of~$\zeta$,
or there exist some intervals of $\zeta$ where it does not hold.
This amounts to study whether the dominant harmonics
can always be found among the primary resonances ($R_1=\bR_1$)
or, on the contrary, secondary resonances have to be taken into account
(and in this case the function $F_1(\zeta)$ is somewhat more complicated).
Such two possiblities also take place in the quadratic case considered
in \cite{DelshamsGG16}.

We can provide an alternative definition for $F_1(\zeta)$
as the minimum of the following functions,
associated to any given resonant sequence $s(q,n)$:
\beq\label{eq:F1q0}
  \wtl F_1^{(q)}(\zeta):=\min_{n\ge0}f^*_{s(q,n)}(\zeta)
\eeq
(for the primary resonances, we have $\wtl F_1^{(\wh q)}=\bF_1$).
Clearly, it is enough to consider essential primitives ($q\in\Pc_0$),
and hence we can write
\beq\label{eq:F1q}
  F_1(\zeta)=\min_{q\in\Pc_0}\wtl F_1^{(q)}(\zeta).
\eeq
Such functions $\wtl F_1^{(q)}(\zeta)$ are completely analogous
to $\bF_1(\zeta)$.
We are going to study only the function $\bF_1(\zeta)$,
showing that it is quasiperiodic and providing lower and upper bounds for it,
and the same will hold for $\wtl F_1^{(q)}(\zeta)$,
with the bounds multiplied by the factor $(\tl\gamma^*_q)^{1/3}\ge1$
in view of~(\ref{eq:fqn1}).
Notice also that only a finite number of primitives $q$ are involved
in~(\ref{eq:F1q}), due to the fact that the (normalized) limits $\tl\gamma^*_q$
have the lower bound~(\ref{eq:lowerbound}),
which is increasing with respect to $\abs q$.

\bremark
Although we implicitly assume that there exists only one sequence
of primary resonances
(see remark~\ref{rk:oneprimary}
at the end of Section~\ref{sect:resonantseq2}),
it is not hard to adapt our definitions and results to the case of
two or more sequences of primary resonances.
In this case, we would choose in~(\ref{eq:s0n}) one of such sequences
as ``the'' sequence $s_0(n)$, when the functions
$\bg_n(\eps)$ and $\bbf_n(\zeta)$
are defined in~(\ref{eq:gn1}) and~(\ref{eq:fqn3})
(see also \cite{DelshamsGG16}).
\eremark

Now we proceed to study the function $\bF_1(\zeta)$ introduced
in~(\ref{eq:bF1}).
Notice that we can regard this function as an $\Ord(\delta)$-perturbation
of the function obtained if we had $\delta=0$ in~(\ref{eq:delta})
(and hence $\bb_n=1$ in~(\ref{eq:gn3})).
Of course, this is fictitious since $\delta$ is determined
by the frequency vector $\omega$ and is not a true parameter.
With this in mind, we define ``unperturbed'' functions
\bea
  \nonumber
  &&\bbf^{(0)}_n(\zeta):=\Cc(\zeta\,;\,n,1)=\Cc_0(\zeta-n),
\\
  \label{eq:F01}
  &&\bF^{(0)}_1(\zeta)
  :=\min_n\bbf^{(0)}_n(\zeta)=\bbf^{(0)}_{N^{(0)}_1}(\zeta).
\eea
The index $N^{(0)}_1=N^{(0)}_1(\zeta)$ providing the minimum
can easily be determined.
On one hand, we use that each function $\bbf^{(0)}_n(\zeta)$
reaches it minimum at $\zeta_n=n$.
On the other hand, applying Lemma~\ref{lm:ZW}(a) (with $Z=1$ and $W=1$)
we find its corners, given by the (transverse) intersection
between the graphs of consecutive functions
$\bbf^{(0)}_n(\zeta)$ and $\bbf^{(0)}_{n+1}(\zeta)$:
\beq\label{eq:xi0}
  \zeta'_n:=n+\xi_0,
  \qquad
  \xi_0:=2\Lg\frac{2\lambda}{\sqrt\lambda+1}\,,
  \quad \textrm{i.e.} \quad
  \lambda^{3\xi_0/2}=\frac{2\lambda}{\sqrt\lambda+1}\,.
\eeq
Hence, we can write $\xi_0=\xi_0(\omega)$ and,
using that $\lambda>1$, it is not hard to see that $1/3<\xi_0<1/2$
(see in Section~\ref{sect:h1_cubicgolden} the concrete value
for the case of the cubic golden vector).
Introducing the intervals $\I_n:=[\zeta'_{n-1},\zeta'_n]$,
we see that $N^{(0)}_1(\zeta)=n$ for any $\zeta\in\I_n$
(strictly speaking, there are two possible values
at the endpoints $\zeta'_n$ of the intervals).
In this way, the function $N^{(0)}_1(\zeta)$ is ``piecewise-constant''
with jump discontinuities at the points $\zeta'_n$,
and the function $\bF^{(0)}_1(\zeta)$ is 1-periodic, continuous
and piecewise-smooth with corners at the same points $\zeta'_n$.
We also obtain the following extreme values:
\bea
  \label{eq:J0unpert}
  \lefteqn{\min\bF^{(0)}_1(\zeta)=\bF^{(0)}_1(n)=\Cc_0(0)\ =\ 1\,,}
\\
  \nonumber
  \lefteqn{\max\bF^{(0)}_1(\zeta)=\bF^{(0)}_1(\zeta'_n)
    =\Cc_0(\xi_0)=\Cc_0(\xi_0-1)}
\\
  \label{eq:J1unpert}
  &&=\ J^{(0)}_1=J^{(0)}_1(\omega)
  \ :=\ \frac13\pq{2\p{\frac{\sqrt\lambda+1}{2\lambda}}^{1/3}
                   +\p{\frac{2\lambda}{\sqrt\lambda+1}}^{2/3}}.
\eea

Returning to the ``perturbed'' function $\bF_1(\zeta)$,
the next lemma shows that, for any $\zeta$,
the index $N_1(\zeta)$ providing the minimum in definition~(\ref{eq:bF1}),
can be found among a finite number (not depending on $\zeta$)
of values around $N^{(0)}_1(\zeta)$.

\begin{lemma}\label{lm:intervalN}
For any $\zeta$, we have
$N^{(0)}_1(\zeta)-N^-\le N_1(\zeta)\le N^{(0)}_1(\zeta)+N^+$,
where we define
\bean
  &&N^-=N^-(\omega):=\log_\lambda\pq{\max\p{
    \frac{1+\delta}{1-\delta}\;,\,2(1+\delta)^{1/2}\lambda^{3(1-\xi_0)/2}+1
  }},
\\
  &&N^+=N^+(\omega):=\log_\lambda\pq{\max\p{
    \frac{1+\delta}{1-\delta}\;,
    \,\p{\frac{\lambda^{3\xi_0/2}+2(1+\delta)^{1/2}}{2(1-\delta)^{1/2}}}^2
  \;}}.
\eean
\end{lemma}

\proof
Let us assume that $\zeta$ belongs to a concrete interval $\I_n$,
where we have $N^{(0)}_1(\zeta)=n$.
In order to show that $N_1(\zeta)$ belongs to the interval $[n-N^-,n+N^+]$,
we have to show that, for any $m$ not belonging to this interval, we~have
\beq\label{eq:fnm1}
  \bbf_m(\zeta)>\bbf_n(\zeta)
  \qquad \textrm{for any} \quad \zeta\in\I_n.
\eeq
To study the relative position of the functions
$\bbf_n(\zeta)$ and $\bbf_m(\zeta)$ (defined in~(\ref{eq:fqn3})),
we will apply Lemma~\ref{lm:ZW} showing that their graphs
do intersect at a point $\zeta^*_{n,m}$, which satisfies:
\beq\label{eq:fnm2}
  \begin{array}{ll}
    \zeta^*_{n,m}<\zeta'_{n-1} &\qquad \textrm{if} \quad m-n<-N^-,
  \\[4pt]
    \zeta^*_{n,m}>\zeta'_n     &\qquad \textrm{if} \quad m-n>N^+,
  \end{array}
\eeq
which says that the (unique) intersection takes place
outside the interval $\I_n$,
and implies the inequality~(\ref{eq:fnm1}).

In order to apply Lemma~\ref{lm:ZW}, we consider the values
$Z=\bzeta_m-\bzeta_n$ and $W=(\bb_m/\bb_n)^{1/3}$,
which satisfy the equality
\beq\label{eq:lambdaZ}
  \lambda^Z=W\,\lambda^{m-n}.
\eeq
On the other hand, recalling that $\bb_n,\bb_m\in[1-\delta,1+\delta]$, we have
$\ds\p{\frac{1-\delta}{1+\delta}}^{1/3}
  \le W
  \le\p{\frac{1+\delta}{1-\delta}}^{1/3}$.

To prove the first assertion of~(\ref{eq:fnm2}),
we use the first bound of Lemma~\ref{lm:ZW}(b), which reads
\[
  \zeta^*_{n,m}<\bzeta_n+2\Lg\frac{2\lambda^{m-n}}{1-\lambda^{m-n}}
  \qquad \textrm{if} \quad \lambda^{m-n}<\min(1,W^{-3}),
\]
where we the equality~(\ref{eq:lambdaZ}) has been taken into account.
By the definition of $N^-$, it is clear that
$\lambda^{n-m}>\dfrac{1+\delta}{1-\delta}\ge\p{\min(1,W^{-3})}^{-1}$.
Moreover, the inequality $\zeta^*_{n,m}<\zeta'_{n-1}$ holds provided
\[
  \Lg\bb_n+2\Lg\frac{2\lambda^{m-n}}{1-\lambda^{m-n}}\le-1+\xi_0.
\]
Replacing $\bb_n$ by $1+\delta$,
the subsequent inequality can be rewritten as
\[
  \lambda^{n-m}\ge2(1+\delta)^{1/2}\lambda^{3(1-\xi_0)/2}+1,
\]
also included in the definition of $N^-$,
which completes the proof of the first assertion of~(\ref{eq:fnm2}).

For the second assertion of~(\ref{eq:fnm2}) we can proceed in similar terms,
using the second bound of Lemma~\ref{lm:ZW}(b).
Nevertheless, the associated computations are somewhat different
due to the lack of symmetry of the functions $\bbf_n(\zeta)$ in the cubic case
(see the remark after the definitions~(\ref{eq:fqn1}--\ref{eq:fqn3})).
We omit the details.
\phantom{\ref{eq:fqn2}}
\qed

In the following proposition, we provide a lower and an upper bound
for the functions $\bF_1(\zeta)$ and $F_1(\zeta)$, and hence for $h_1(\eps)$,
as $\Ord(\delta)$-perturbations of the values obtained
in~(\ref{eq:J0unpert}--\ref{eq:J1unpert}).
More precisely, such bounds will be given by the values
\beq\label{eq:J01}
  J^-_0=J^-_0(\omega):=(1-\delta)^{1/3},
  \qquad
  J^+_1=J^+_1(\omega):=J^{(0)}_1\,(1+\delta)^{1/3},
\eeq
which satisfy $0<J^-_0<1<J^{(0)}_1<J^+_1$.
Recall that lower and an upper bounds for $h_1(\eps)$ or,
equivalently, for $F_1(\zeta)$,
can be associated to upper and lower bounds
for the splitting distance, respectively
(see also \cite{DelshamsGG14a}).
Recalling the value $B^-_0=B^-_0(\omega)$ defined in~(\ref{eq:Bm0}),
we also introduce the ``\emph{strong separation condition}'':
\beq\label{eq:strongsep}
  B^-_0\ \ge\ J^+_1,
\eeq
which is somewhat more restrictive than
the ``weak separation condition'' introduced in~(\ref{eq:weaksep}).
Under the strong condition, the inequality~(\ref{eq:FF1}) becomes an equality,
i.e.~the dominant harmonic is always given by a primary resonance,
and hence the function $F_1(\zeta)=h_1(\eps)$ becomes somewhat simpler.
Such a condition is fulfilled for the cubic golden frequency vector,
as we show in Section~\ref{sect:h1_cubicgolden}.

\begin{proposition}\label{prop:J01}
The functions $F_1(\zeta)$ and $\bF_1(\zeta)$ are positive, continuous and
piecewise-smooth, and satisfy for any $\zeta$ the bounds:
\[
  J^-_0\le F_1(\zeta)\le\bF_1(\zeta)\le J^+_1,
\]
with $J^-_0$ and $J^+_1$ defined in~(\ref{eq:J01}).
Moreover, if the strong separation condition~(\ref{eq:strongsep}) is fulfilled,
then we have $F_1(\zeta)=\bF_1(\zeta)$ for any $\zeta$, and hence
the most dominant harmonic is always given by a primary resonance.
\end{proposition}

\proof
The lower bound for $F_1(\zeta)$ is a direct consequence
of~(\ref{eq:F1q0}--\ref{eq:F1q}), using that for any $k=s(q,n)\in\A$
we have the lower bound
\beq\label{eq:fsqn_lower}
  f^*_{s(q,n)}(\zeta)\ge(\tl\gamma^*_qb_{s(q,n)})^{1/3}\ge(1-\delta)^{1/3},
\eeq
which comes from~(\ref{eq:fqn1}), using also that
$b_{s(q,n)}\ge1-\delta$ by~(\ref{eq:bsqn}).

To provide an upper bound for $\bF_1(\zeta)$,
we take into account that $\bb_n\le1+\delta$ and introduce the function
\[
  \bF^{\,+}_1(\zeta):=\min_{n\ge0}\bbf^+_n(\zeta),
  \qquad
  \bbf^+_n(\zeta):=\Cc(\zeta\,;\,\bzeta^+_n,1+\delta),
  \quad
  \bzeta^+_n:=n+\Lg(1+\delta),
\]
defined as in~(\ref{eq:bF1}) but replacing $\bb_n$ by $1+\delta$
in~(\ref{eq:fqn3}).
Notice that the function $\bF^{\,+}_1(\zeta)$
can easily be related to the ``unperturbed'' function defined
in~(\ref{eq:F01}): for any $\zeta$, we have
\[
  \bF^{\,+}_1(\zeta)=(1+\delta)^{1/3}\,\bF^{(0)}_1(\zeta-\Lg(1+\delta)),
\]
and we deduce from~(\ref{eq:J1unpert}) and~(\ref{eq:J01})
that $\max\bF^{\,+}_1(\zeta)=J^+_1$.

We study the relative position of the graphs of the functions
$\bbf_n(\zeta)$ and $\bbf^+_n(\zeta)$ by applying Lemma~\ref{lm:ZW}(a), with
$Z=\bzeta^+_n-\bzeta_n=\Lg((1+\delta)/\bb_n)$ and $W=((1+\delta)/\bb_n)^{1/3}$.
In general we have $\bb_n<1+\delta$ and,
since $\lambda^Z=W$, the graphs do not intersect and we have
$\bbf_n(\zeta)<\bbf^+_n(\zeta)$ for any $\zeta$.
Instead, if $\bb_n=1+\delta$ (a rather particular case)
then the two functions obviously coincide.
We deduce, for any $\zeta$, the bound
\beq\label{eq:F1_upper}
  \bF_1(\zeta)\le\bF^{\,+}_1(\zeta)\le J^+_1\,.
\eeq

Finally, to show that the strong separation condition~(\ref{eq:strongsep})
implies the equality $F_1(\zeta)=\bF_1(\zeta)$,
it is enough to see that a lower bound for the functions $\wtl F_1^{(q)}(\zeta)$
introduced in~(\ref{eq:F1q0}), for $q\ne\wh q$,
is greater than the upper bound $J^+_1$ for $\bF_1(\zeta)$, obtained above.
Indeed, for secondary resonances $s(q,n)$, with $q\ne\wh q$,
the lower bound~(\ref{eq:fsqn_lower}) becomes 
\[
  f^*_{s(q,n)}(\zeta)\ge\p{\tl\gamma^*_{\whh q}(1-\delta)}^{1/3}=B^-_0\ge J^+_1\,,
\]
where $\gamma^*_{\whh q}$ is the minimum of the ``mean Diophantine constants''
for secondary resonances (see~(\ref{eq:gammasec})),
and the same lower bound holds
for the functions $\wtl F_1^{(q)}(\zeta)$, $q\ne\wh q$.
\qed

\bremark
It is an interesting question whether the lower and upper bounds
$J^-_0$ and $J^+_1$ provided by this proposition are \emph{sharp},
i.e.~they coincide with the infimum and the supremum
of the function $F_1(\zeta)$.
On one hand, we can expect the lower bound $J^-_0$
(and hence the upper bound for the splitting)
to be sharp, since for primary resonances the
lower bounds~(\ref{eq:fsqn_lower}) are given by the factors $\bb_n$,
which will can be arbitrarily close to $1-\delta$ for suitable $n$.
Instead, in general the upper bound $J^+_1$
(and hence the lower bound for the splitting)
is far from being sharp, because it has been obtained in~(\ref{eq:F1_upper})
by considering, for all $n$,
the worst possible case in the bound $\bb_n\le1+\delta$.
In Section~\ref{sect:qp}, we prove the sharpness of the lower bound $J^-_0$
and show that, for a given frequency vector $\omega$,
we can give (numerically) a~sharp upper bound $J^*_1$ ($\le J^+_1$),
using the quasiperiodicity of the function $F_1(\zeta)$.
In the same way, it would be enough to assume that
$B^-_0\ge J^*_1$\,, instead of~(\ref{eq:strongsep}),
in order to ensure that the splitting can be described
in terms of only the primary resonances.
This value $J^*_1$ is computed in Section~\ref{sect:h1_cubicgolden}
for the concrete case of the cubic golden frequency vector.
\eremark

To end this section, we also deduce some useful
properties of the function $S_1=S_1(\eps)$, giving the dominant harmonic.
Namely, this function is ``piecewise-constant'',
with jump discontinuities exactly at the corners of $h_1(\eps)$.
Moreover, its asymptotic behavior as $\eps\to0$ turns out to be polynomial:
\beq\label{eq:estimS}
  \abs{S_1(\eps)} \sim \frac1{\eps^{1/6}}\,.
\eeq
Indeed, the most dominant harmonic belongs to some resonant sequence:
we can write $S_1(\eps)=s(q,N)$
for some $q=q(\eps)$, and for $N=N(\eps)$ such that the value $\eps^*_{s(q,N)}$
is close to $\eps$, among the sequence $\eps^*_{s(q,n)}$, $n\ge0$.
Recalling~(\ref{eq:gqn2}) and the estimate
$\abs{s(q,N)}\sim\lambda^{N/2}=\p{\lambda^{3N}}^{1/6}$
deduced from~(\ref{eq:cubicfreq_a}), we get~(\ref{eq:estimS}).
Notice that it is not necessary to include $q$
in the estimate~(\ref{eq:estimS})
(in spite of the fact that $K_q$ and $\tl\gamma^*_q$
appear in the expression~(\ref{eq:gqn2})),
since only a finite number of resonant sequences $s(q,\cdot)$ is involved.

\subsection{Quasiperiodicity of the estimate of the most
  dominant harmonic}
\label{sect:qp}

Now, our aim is to show that the function $F_1(\zeta)$
is \emph{quasiperiodic} with frequencies 1 and $\phi$.
As we show below, this property is directly related
to the oscillating factors $b_{s(q,n)}$ introduced in~(\ref{eq:bsqn})
for each resonant sequence,
denoted~$\bb_n$ in~(\ref{eq:gn3})
for the particular case of the primary resonances.
Moreover, the facts that $\phi$ is an irrational number
by Lemma~\ref{lm:irrational}, and $\delta>0$ by Lemma~\ref{lm:delta},
allow us to ensure that  the function $F_1(\zeta)$ is \emph{not periodic},
which makes an important difference
with respect to the case of quadratic frequencies
considered in \cite{DelshamsGG16}.

Recall that, in~(\ref{eq:F1q}), we wrote $F_1(\zeta)$ as the minimum
of the functions $\wtl F_1^{(q)}(\zeta)$,
associated to each resonant sequence~$s(q,n)$.
Since all such functions are analogous to the function $\bF_1(\zeta)$,
associated to the primary resonances~$s_0(n)$ and defined in~(\ref{eq:bF1}),
it is enough to show the quasiperiodicity of $\bF_1(\zeta)$.

As a rough explanation for the frequencies 1 and $\phi$,
notice that we can consider $\bF_1(\zeta)$ as an $\Ord(\delta)$-perturbation
of the function $\bF^{(0)}_1(\zeta)$ introduced in~(\ref{eq:F01}),
which is 1-periodic with respect to $\zeta$,
and the oscillating factors~$\bb_n$ defined in~(\ref{eq:gn3})
give rise to the second frequency~$\phi$.

\begin{figure}[!b]
  \centering
  \includegraphics[width=0.4\textwidth,angle=-90]{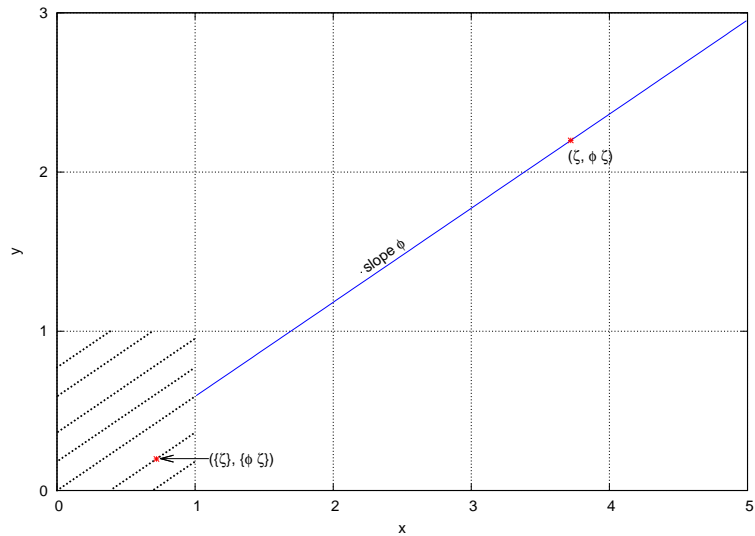}
  \caption{\small\emph{
    The function $\Upsilon(x,y)$ on $\R^2$ ``interpolating'' $\bF_1(\zeta)$
    along the straight lines $x=\zeta$, $y=\phi\,\zeta$,
    and its reduction to the torus~$\T_*^{\,2}$
    (the slope $\phi\approx0.590935$ corresponds to the case
    of the cubic golden vector).}}
\label{fig:interp}
\end{figure}

To be more precise, we are going to construct a positive, continuous
and piecewise-smooth function
$\Upsilon(x,y)$, defined on $\R^2$ and 1-periodic with respect to $x$ and $y$,
such that
\beq\label{eq:qp1}
  \Upsilon(\zeta,\phi\,\zeta)=\bF_1(\zeta)
  \qquad \textrm{for any} \quad \zeta\ge\zeta_0
\eeq
(for some $\zeta_0$ to be determined below, in Proposition~\ref{prop:upsilon}).
Equivalently, we can consider $\Upsilon(x,y)$ as defined on
a torus~$\T_*^{\,2}$, with $\T_*:=\R/\Z$ represented as the interval $[0,1)$,
and the above equality can be rewritten as
\bea
  \label{eq:qp2}
  \lefteqn{\Upsilon(\zeta,\pp{\phi(j+\zeta)})=\bF_1(j+\zeta)}
\\
  \nonumber
  &&\textrm{for any integer} \ j\ge0 \ \textrm{and} \ \zeta\in[0,1),
  \ \textrm{with} \ j+\zeta\ge\zeta_0.
\eea
where $\pp a\in[0,1)$ denotes the fractional part of a given number $a\in\R$.
This property of ``interpolation'' is illustrated in Figure~\ref{fig:interp}.

Like $\bF_1(\zeta)$, defined in~(\ref{eq:bF1}) as the minimum of the
functions $\bbf_n(\zeta)$, the ``interpolating'' function $\Upsilon(x,y)$
will be defined in a similar way, as the minimum of a family functions.
First of all, we define the 1-periodic function
\[
  \beta(y):=1+\delta\cos(2\pi\cdot y+2\psi_{\wh q}-\theta), \qquad y\in\R,
\]
and it is clear that the oscillating factors~(\ref{eq:gn3})
are ``interpolated'' by this function: $\beta(\pp{n\phi})=\bb_n$ for any $n$
(we can say that the values $\pp{n\phi}$, filling densely the circle $\T_*$,
are replaced by the continuous variable $y$).
Now, recalling the ``hyperbolic cosine-like'' functions $\Cc(\zeta\,;\,Z,Y)$
introduced in~(\ref{eq:defC}), we define for $n\in\Z$ the functions
\beq\label{eq:chin}
  \chi_n(x,y)
  :=\Cc(x\,;\,n+\Lg\beta(y-\phi\,x+\pp{n\phi})\,,
      \,\beta(y-\phi\,x+\pp{n\phi})),
   \qquad (x,y)\in\R^2,
\eeq
which are clearly smooth and 1-periodic with respect to $y$,
but not periodic with respect to $x$.
Finally, we define
\beq\label{eq:upsilon}
  \Upsilon(x,y):=\min_{n\in\Z}\chi_n(x,y)=\chi_{\wtl N_1}(x,y),
  \qquad (x,y)\in\R^2,
\eeq
with $\wtl N_1=\wtl N_1(x,y)$ (compare with~(\ref{eq:bF1})).

It is clear that the functions $\chi_n(x,y)$ are closely related
to the functions $\bbf_n(\zeta)$ defined in~(\ref{eq:fqn3}),
as we see from the definition~(\ref{eq:chin}),
by restricting $(x,y)$ to straight lines of slope $\phi$.
To express this relationship more clearly we define, for any $y_0\in\R$,
a function of one variable by restricting $\chi_n(x,y)$
to any straight line $y=y_0+\phi\,x$ for a given $y_0$,
\bea
  \label{eq:chin1}
  &&\wh\chi_n(x\,;\,y_0)
  :=\chi_n(x,y_0+\phi\,x)=\Cc(x\,;\,\bx_n(y_0),\bbeta_n(y_0)),
\\
  \nonumber
  &&\bx_n(y_0):=n+\Lg\bbeta_n(y_0),
  \qquad
  \bbeta_n(y_0):=\beta(y_0+\pp{n\phi})
\eea
(compare with~(\ref{eq:fqn3})).
We can also define
\beq\label{eq:upsilon1}
  \wh\Upsilon(x\,;\,y_0)
  :=\min_{n\in\Z}\wh\chi_n(x\,;\,y_0)
  =\wh\chi_{\wh N_1}(x\,;\,y_0),
\eeq
and it is clear that
$\wh\Upsilon(x\,;\,y_0)=\Upsilon(x,y_0+\phi\,x)$,
and also $\wh N_1(x\,;\,y_0)=\wtl N_1(x,y_0+\phi\,x)$
(with the difference that $\Upsilon$ is 1-periodic
and can be reduced to $\T_*^{\,2}$, see Proposition~\ref{prop:upsilon},
but the periodicity with respect to $x$
does not hold for $\wh\Upsilon$).

Some of the properties stated in the following lemma
are clearly inherited from the results
of Lemmas~\ref{lm:fkk}, \ref{lm:ZW} and~\ref{lm:intervalN}.

\begin{lemma}\label{lm:chin}
\ \btm
\item[(a)] The functions $\chi_n(x,y)$ are smooth and 1-periodic
with respect to $y$, and satisfy the following translation property:
\[
  \chi_n(x+1,y)=\chi_{n-1}(x,y),
  \qquad \textrm{for any} \quad (x,y)\in\R^2, \quad n\in\Z.
\]
\item[(b)] For any given $n$ and $y_0\in\R$,
the function $\wh\chi_n(x\,;\,y_0)$
is convex (with respect to $x$) and attains its minimum at $x=\bx_n(y_0)$,
with the minimum value $\bbeta_n(y_0)^{1/3}$.
The dependence of $\wh\chi_n(x\,;\,y_0)$
on the parameter $y_0$ is~1-periodic.
\item[(c)]  For any given $n$, the function $\chi_n(x,y)$
attains its minimum at the point $(x,y)=(\tl x_n,\tl y_n)$, with
\[
  \tl x_n=n+\Lg(1-\delta),
  \quad
  \tl y_n\equiv\frac{\pi-2\psi_{\wh q}+\theta}{2\pi}+\phi\,\Lg(1-\delta) \pmod1,
\]
with the minimum value $(1-\delta)^{1/3}$.
\item[(d)] For any given $n,m$ with $n\ne m$, and $y_0\in\R$,
the functions $\wh\chi_n(x\,;\,y_0)$ and $\wh\chi_m(x\,;\,y_0)$
do not coincide. Their graphs intersect transversely at a unique point,
or do not intersect.
The set $\Y_{n,m}$ of values $y_0$ such that the intersection
exists is a union of open intervals
(or eventually $\Y_{n,m}=\R$, $\Y_{n,m}=\emptyset$).
For $y_0\in\Y_{n,m}$, the intersecting point $x=x^*_{n,m}(y_0)$
(given explicitly in~(\ref{eq:xnm}))
is a smooth and 1-periodic function of $y_0$.
\item[(e)] For any given $n,m$ with $n\ne m$, the graphs of the functions
$\chi_n(x,y)$ and $\chi_m(x,y)$ intersect (if they do) transversely
along the curves parameterized by
\[
  x=x^*_{n,m}(y_0), \quad y=y_0+\phi\,x^*_{n,m}(y_0), \qquad y_0\in\Y_{n,m}.
\]
\item[(f)]
For any $(x,y)$, we have
$N^{(0)}_1(x)-N^-\le\wtl N_1(x,y)\le N^{(0)}_1(x)+N^+$,
with $N^{(0)}_1(x)$ as in~(\ref{eq:F01}),
and $N^\pm=N^\pm(\omega)$ as in Lemma~\ref{lm:intervalN}.
\etm
\end{lemma}

\proof
The only assertion to be checked in~(a) is the translation property.
For that, it is enough to ensure that
\[
  \beta(y-\phi(x+1)+\pp{n\phi})=\beta(y-\phi\,x+\pp{(n-1)\phi}),
\]
but this is a direct consequence of the 1-periodicity of $\beta(y)$.
The proof of~(b) is straightforward from the definition of the functions
$\wh\chi_n(x\,;\,y_0)$ in~(\ref{eq:chin1}).
We also get~(c) as a direct consequence of~(b),
choosing $y_0=y^{(n)}_0$ such that $\bbeta_n(y_0)$
attains its minimum value $1-\delta$, and hence $\tl x_n=\bx_n(y^{(n)}_0)$,
$\tl y_n=y^{(n)}_0+\phi\,\tl x_n$.

For~(d), we first notice that the functions $\wh\chi_n(x\,;\,y_0)$
and $\wh\chi_m(x\,;\,y_0)$ do not coincide,
since $\bbeta_n(y_0)\ne\bbeta_{m}(y_0)$ (due~to the irrationality of~$\phi$).
Then, we directly apply Lemma~\ref{lm:ZW} with
$Z=\bx_{m}(y_0)-\bx_n(y_0)$ and $W=(\bbeta_{m}(y_0)/\bbeta_n(y_0))^{1/3}$.
We get the formula for the intersecting point,
\beq\label{eq:xnm}
  x^*_{n,m}(y_0)
  =\bx_n(y_0)+2\Lg\frac{2\lambda^Z(W\lambda^{Z/2}-1)}{\lambda^Z-W}\,.
\eeq
If the intersection exists, it is unique, but its existence
may depend on $y_0$, according to the condition given in Lemma~\ref{lm:ZW}.
We also get~(e) as a direct consequence of~(d).

Finally, for the proof of~(f), for any $y_0$ we consider the function
$\wh\Upsilon(x\,;\,y_0)$ defined in~(\ref{eq:upsilon1}),
and it is enough to prove that
$N^{(0)}_1(x)-N^-\le\wh N_1(x\,;\,y_0)\le N^{(0)}_1(x)+N^+$.
Now, we can use that the functions $\wh\chi_n(x\,;\,y_0)$
introduced in~(\ref{eq:chin1})
are completely analogous to the functions $\bbf_n(\zeta)$ in~(\ref{eq:fqn3}),
replacing $\bb_n$ by $\bbeta_n(y_0)$, and $\bzeta_n$ by $\bx_n(y_0)$.
Then, the proof follows exactly as in Lemma~\ref{lm:intervalN},
using the values of $Z$ and $W$ defined above.
\qed

\begin{proposition}\label{prop:upsilon}
The function $\Upsilon(x,y)$ is continuous and piecewise-smooth,
and 1-periodic with respect to $x$ and $y$,
and satisfies the ``interpolation'' property~(\ref{eq:qp1})
for $\zeta\ge\zeta_0:=N^-+\xi_0$
(recall that $\xi_0$ is defined in~(\ref{eq:xi0})).
\end{proposition}

\proof
First of all, from definitions~(\ref{eq:fqn3}) and~(\ref{eq:chin}),
it is not hard to see that the equality
$\chi_n(\zeta,\phi\,\zeta)=\wh\chi_n(\zeta\,;\,0)=\bbf_n(\zeta)$
is fulfilled for any $n\ge0$ and $\zeta\in\R$
(we only have to use that $\bbeta_n(0)=\bb_n$).
By Lemma~\ref{lm:chin}(f),
we can take the minimum over $n$ by restricting ourselves
to a finite number of cases,
$N^{(0)}_1(\zeta)-N^-\le n\le N^{(0)}_1(\zeta)+N^+$,
and we directly get the equality~(\ref{eq:qp1}),
or equivalently~(\ref{eq:qp2}).
However, in order to ensure that $n\ge0$ as in the definition~(\ref{eq:bF1}),
we need that $N^{(0)}_1(\zeta)\ge N^-$.
As can be seen in~(\ref{eq:F01}), we have $N^{(0)}_1(\zeta)\ge\zeta-\xi_0$,
and hence we assume $\zeta\ge N^-+\xi_0$.

The fact that $\Upsilon(x,y)$ is, for any $(x,y)$,
the minimum of a finite number of smooth functions ensures
that it is continuous and piecewise-smooth.
It is also clear that it is periodic with respect to $y$,
since so are the functions $\chi_n(x,y)$.
Finally, its periodicity with respect to $x$
is easily deduced from the translation property of Lemma~\ref{lm:chin}(a).
\qed

\begin{figure}[!b]
  \centering
  \subfigure{
    \includegraphics[width=0.37\textwidth,angle=-90]{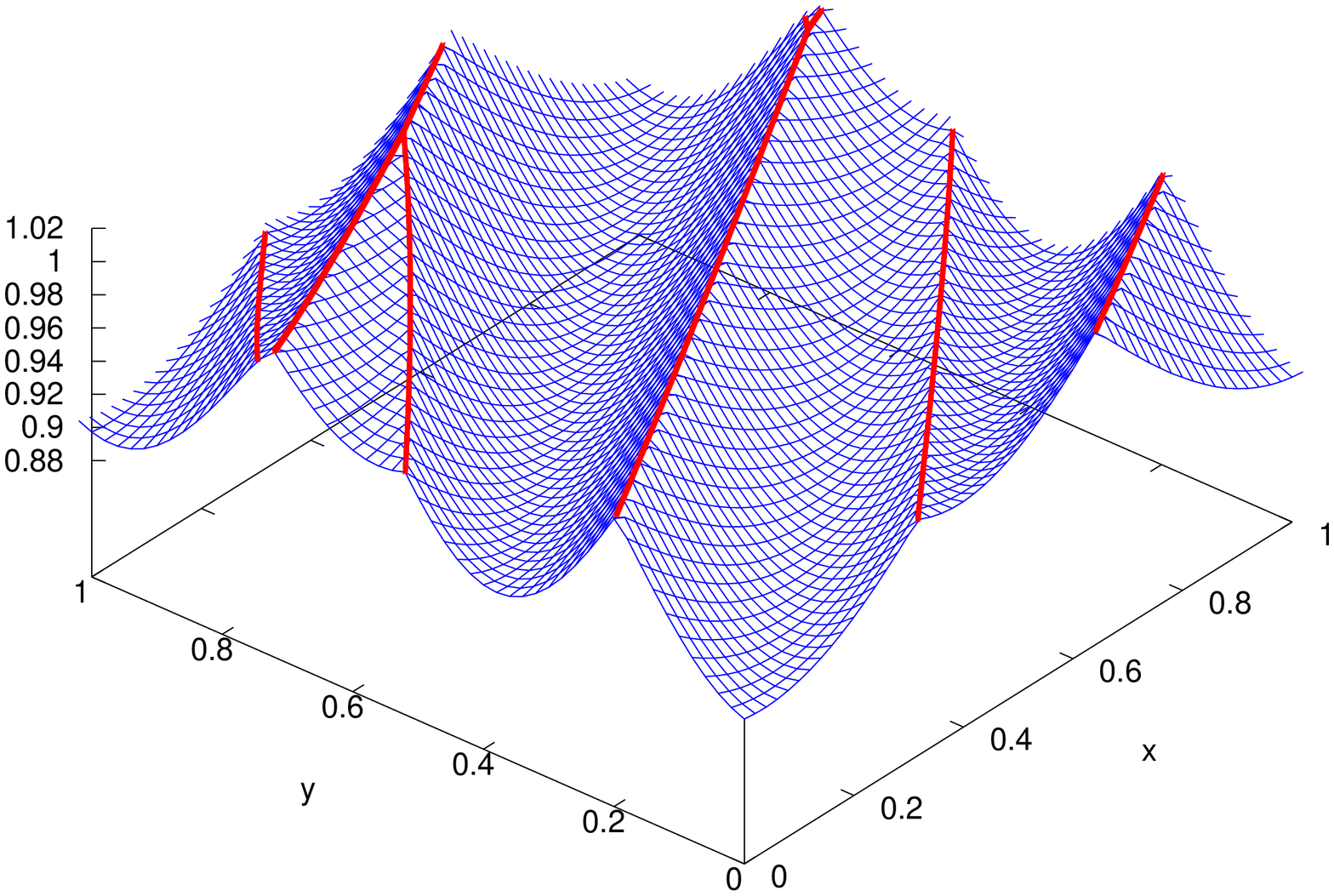}}
  \hspace{-1.2cm}
  \subfigure{
    \includegraphics[width=0.35\textwidth,angle=-90]{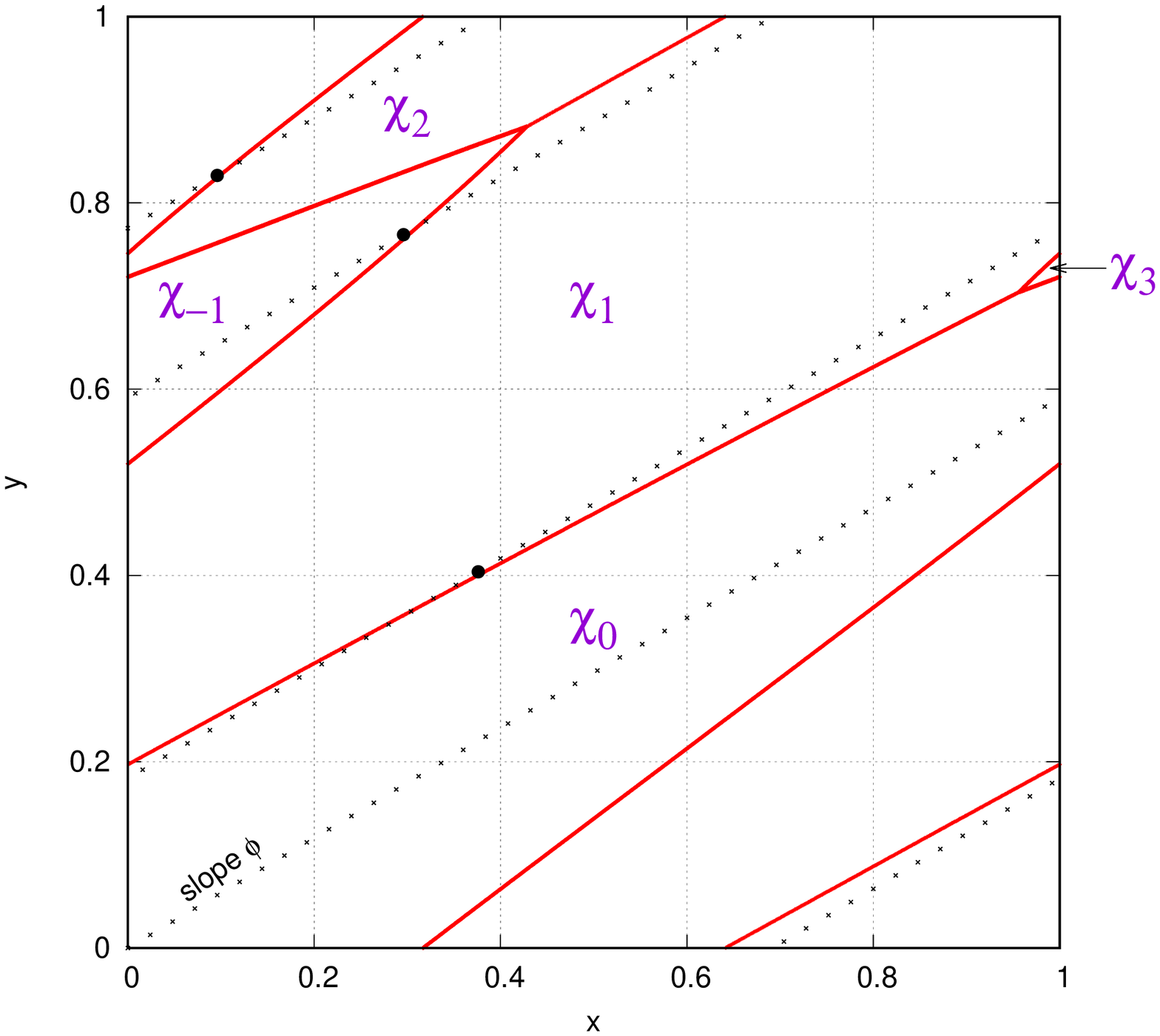}}
  \caption{\small\emph{
    Graph of the function $\Upsilon(x,y)$ on its domain $\T_*^{\,2}$,
    as the minimum of the functions $\chi_n(x,y)$, for the cubic golden vector.
    The red curves (which are not straight lines) are the borders between the
    regions of dominance, where a different function $\chi_n(x,y)$ gives the
    minimum in~(\ref{eq:upsilon}). The function $\bF_1(\zeta)$ is the
    restriction of $\Upsilon(x,y)$ along the dotted line of slope $\phi$, by
    the property of ``interpolation'', see~(\ref{eq:qp1}--\ref{eq:qp2}). The
    changes in the dominance, which take place when the line of slope $\phi$
    crosses a red curve, correspond to the corners of $\bF_1(\zeta)$ in
    Figure~\ref{fig:h1}.}}
\label{fig:upsilon}
\end{figure}

In this way, by studying the function $\Upsilon(x,y)$ on the torus $\T_*^{\,2}$
we can determine the intervals of dominance for the function $\bF_1(\zeta)$,
in~(\ref{eq:bF1}).
It is enough to divide $\T_*^{\,2}$ into a finite number of regions, according
to the function $\chi_n(x,y)$ giving the minimum in~(\ref{eq:upsilon}).
Since for $x\in[0,1)$ the index $N^{(0)}_1(x)$ is either 0 or 1,
by Lemma~\ref{lm:chin}(f) it is enough to consider
the functions $\chi_n(x,y)$ with $-N^-\le n\le1+N^+$.
The regions visited by the straight line $(\zeta,\phi\,\zeta)$
correspond the intervals of dominance for $\bF_1(\zeta)$.
See Figure~\ref{fig:upsilon} for an illustration, for the concrete case
of the cubic golden vector
(we point out that the borders between neighbor regions are not straight lines,
but rather pieces of the curves parameterized in Lemma~\ref{lm:chin}(e)).

Numerically, we can obtain \emph{sharp bounds}
for the function $\bF_1(\zeta)$, improving the ones given
in Proposition~\ref{prop:J01}. Since $\phi$ is irrational,
the line $(\zeta,\phi\,\zeta)$ fills densely the torus $\T_*^{\,2}$
and hence
\[
  \inf\bF_1(\zeta)=\min\Upsilon(x,y)=J^-_0,
  \qquad
  \sup\bF_1(\zeta)=\max\Upsilon(x,y)\le J^+_1.
\]
The minimum value $J^-_0=(1-\delta)^{1/3}$ of $\Upsilon(x,y)$
is attained at the point given in Lemma~\ref{lm:chin}(c),
choosing $n$ such that $\tl x_n\in[0,1)$.
On the other hand, by the convexity of $\Upsilon$ along the lines
of slope $\phi$, the maximum value
\beq\label{eq:J1ast}
  J^*_1:=\max\Upsilon(x,y)
\eeq
is attained at some point belonging to some of the curves
limiting the regions of dominance illustrated in Figure~\ref{fig:upsilon},
Recall that the values $J^-_0$ and $J^*_1$ are associated, respectively,
to sharp upper and lower bounds for the maximum splitting distance
(see remark~\ref{rk:supremum} after Theorem~\ref{thm:main}).
Again, see Section~\ref{sect:h1_cubicgolden}
for the case of the cubic golden vector.

\subsection{The particular case of the cubic golden frequency
  vector}
\label{sect:h1_cubicgolden}

As a continuation of Section~\ref{sect:cubicgolden},
we provide particular data concerning the function $h_1(\eps)=F_1(\zeta)$,
and hence the asymptotic estimate for the splitting,
for the concrete case of the cubic golden frequency vector
introduced in~(\ref{eq:cubicgolden}).

First of all, recall that the function $\bF_1(\zeta)$
defined in~(\ref{eq:bF1}), associated to the primary resonances,
is an $\Ord(\delta)$-perturbation
of the 1-periodic function $\bF^{(0)}_1(\zeta)$ introduced in~(\ref{eq:F01}).
This one reaches its minimum value at the points $\zeta_n=n$,
and its maximum value at the points $\zeta'_n=n+\xi_0$,
with $\xi_0\approx0.492049$ in~(\ref{eq:xi0}),
where we have used the value of $\lambda$
obtained in~(\ref{eq:lambdagolden}).
The minimum value is 1 and the maximum value
is $J^{(0)}_1\approx1.009141$ by~(\ref{eq:J1unpert}).

For the ``perturbed'' function $\bF_1(\zeta)$, 
we use the value of $\delta$ obtained in~(\ref{eq:deltagolden}) and,
in Lemma~\ref{lm:intervalN}, we get the values
$N^-\approx3.65$ and $N^+\approx3.97$\,.
This says that, for $\zeta$ belonging
to a given interval $\I_n=[\zeta'_{n-1},\zeta'_n]$
(where we have $N^{(0)}_1(\zeta)=n$), we can compute $\bF_1(\zeta)$
as the minimum of the functions $\bbf_j(\zeta)$ for $n-3\le j\le n+3$.

On the other hand, by Proposition~\ref{prop:J01} we have the following
lower and upper bounds for $\bF_1(\zeta)$,
\[
  J^-_0\approx0.892341\,,
  \quad
  J^+_1\approx1.098383\,.
\]
The strong separation condition~(\ref{eq:strongsep}) is fulfilled
for the cubic golden vector, since the value $B^-_0$
obtained in~(\ref{eq:JB0golden}) is clearly greater than $J^+_1$\,,
and hence $\bF_1(\zeta)=F_1(\zeta)$ for this example.
In fact, the upper bound $J^+_1$ can be replaced by the sharp upper bound
$J^*_1$ defined in~(\ref{eq:J1ast}), and numerically we see that
\[
  J^*_1\approx1.010619
\]
(this value is reached at the confluence of the regions where
$\chi_{-1}$, $\chi_1$, $\chi_2$ are dominant, see Figure~\ref{fig:upsilon}).

\section{Justification of the asymptotic estimate}
\label{sect:technical}

We consider in this section the final step in the proof of our main result
(Theorem~\ref{thm:main}), which gives an exponentially small asymptotic
estimate for the maximal splitting distance.
This requires to bound the sum
of the non-dominant terms of the Fourier expansion of
the Melnikov potential $L(\theta)$,
ensuring that it can be approximated by its dominant harmonic.
Furthermore, to ensure that the Poincar\'e--Melnikov
method~(\ref{eq:melniapproxM})
predicts correctly the size of the splitting in the singular case $\mu=\eps^r$,
we extend the results to the splitting function $\M(\theta)$
by showing that the asymptotic estimate of the dominant harmonic
is large enough to overcome the harmonics of the error term
in~(\ref{eq:melniapproxM}).
This step is analogous to the case of the quadratic golden number
done in \cite{DelshamsG04} (see also \cite{DelshamsGG16}),
using the upper bounds for the error term provided in \cite{DelshamsGS04},
and we omit many details.
In fact, the specific arithmetic properties of cubic frequency vectors
are not used in this section.

We start with describing our approach in a few words.
First of all, notice that Theorem~\ref{thm:main} is stated in terms
of the splitting function $\M=\nabla\Lc$ introduced in~(\ref{eq:defM}).
We write, for the splitting potential and function,
\beq\label{eq:Lk}
  \Lc(\theta)=\sum_{k\in\Zc\setminus\pp0}\Lc_k\,\cos(\scprod k\theta-\tau_k),
  \qquad
  \M(\theta)=-\!\!\sum_{k\in\Zc\setminus\pp0}\M_k\,\sin(\scprod k\theta-\tau_k),
\eeq
with scalar positive coefficients $\Lc_k$,
and vector coefficients
\beq\label{eq:Mk}
  \M_k=k\,\Lc_k\;\in\R^3.
\eeq
Although the Melnikov approximation~(\ref{eq:melniapproxM})
is in principle valid for real $\theta$, it is standard to see that it
can be extended to a complex strip of suitable width
(see for instance \cite{DelshamsGS04}), from which one gets
upper bounds for $\abs{\Lc_k-\mu L_k}$ and $\abs{\tau_k-\sigma_k}$
(see~(\ref{eq:coefL})),
which imply the estimates given below in Lemma~\ref{lm:dominantsL},
ensuring that the most dominant harmonic of the Melnikov potential $L(\theta)$,
obtained for $k=S_1(\eps)$ (see~(\ref{eq:h1})),
is also the dominant one for the splitting potential~$\Lc(\theta)$.
Then, this dominant harmonic determines the asymptotic estimate for the maximal
splitting distance, given in Theorem~\ref{thm:main}.

With this idea, we consider the approximation of $\Lc(\theta)$ given
by its dominant harmonic, as well as the corresponding remainder,
\bea
  &&\nonumber
  \Lc(\theta)=\Lc^{(1)}(\theta)+\F^{(2)}(\theta),
\\
  &&\label{eq:L1F2}
  \Lc^{(1)}(\theta):=\Lc_{S_1}\cos(\scprod{S_1}\theta-\tau_{S_1}),
  \quad
  \F^{(2)}(\theta):=\sum_{k\in\Zc_2}\Lc_k\cos(\scprod k\theta-\tau_k),
\eea
where we denote
$\Zc_2:=\Zc\setminus\{0,S_1\}$,
and we give below, in Lemma~\ref{lm:dominantsL},
an estimate for the sum of all harmonics in the remainder $\F^{(2)}(\theta)$,
in order to ensure that the maximal splitting distance can be approximated
by the size of the coefficient of the most dominant harmonic $S_1(\eps)$.
In fact, the estimate for $\F^{(2)}(\theta)$ is also given,
by the exponential smallness of the harmonics,
in terms of its own dominant harmonic in the set $\Zc_2$,
that we denote as $S_2(\eps)$.
With this in mind, we introduce as in~(\ref{eq:h1}) the
continuous and piecewise-smooth function
\beq\label{eq:h2}
  h_2(\eps):=\min_{k\in\A\setminus\pp{S_1}}g^*_k(\eps)=g^*_{S_2}(\eps),
\eeq
It is not hard to see from Lemmas~\ref{lm:fkk} and~\ref{lm:ZW} that
the corners of $h_1(\eps)$, at which a change in the first dominant harmonic
takes place, are exactly the points $\ceps$ such that $h_1(\ceps)=h_2(\ceps)$
(such points are also the ``lower corners'' of $h_2(\eps)$,
but this function also has ``upper corners'' where it coincides
with the analogous function $h_3(\eps)$ associated
to the third dominant harmonic; see \cite{DelshamsGG16}).

The following lemma, analogous
to the one established in \cite{DelshamsG03,DelshamsG04},
provides an asymptotic estimate for the dominant harmonic $\Lc_{S_1}$,
and an upper bound for the difference of the phase $\tau_{S_1}$
with respect to the original one~$\sigma_{S_1}$,
as well as an estimate for the sum of all the harmonics
in the remainder appearing in~(\ref{eq:L1F2}),
In fact, we are not directly interested in the
splitting potential $\Lc(\theta)$, but rather its derivative $\M(\theta)$.
Recall that the coefficients $\Lc_k$, introduced in~(\ref{eq:Lk}),
are all positive, and that the constant $C_0$ in the exponentials has been
defined in~(\ref{eq:beta_gk}).
On the other hand, we use the following notation:
for positive quantities, we write $f\preceq g$ if we can bound $f\leq c\,g$
with some (positive) constant $c$ not depending on~$\eps$ and $\mu$.
In this way, we can write $f \sim g$ if $g\preceq f\preceq g$.

\begin{lemma}\label{lm:dominantsL}
For $\eps$ small enough and $\mu=\eps^r$ with $r>3$, one has:
\btm
\item[\rm(a)]
$\ds\Lc_{S_1}
 \sim\mu\,L_{S_1}
 \sim\frac{\mu}{\eps^{1/6}}
   \,\exp\pp{-\frac{C_0h_1(\eps)}{\eps^{1/6}}}$,
\quad
$\abs{\tau_{S_1}-\sigma_{S_1}}\preceq\dfrac\mu{\eps^3}$\,;
\vspace{4pt}
\item[\rm(b)]
$\ds\sum_{k\in\Zc_2}\Lc_k
\sim\frac1{\eps^{1/3}}\,\Lc_{S_2}
\sim\frac\mu{\eps^{1/3}}
  \,\exp\pp{-\frac{C_0 h_2(\eps)}{\eps^{1/6}}}$.
\etm
\end{lemma}

\sketchproof
We only give the main ideas of the proof, since it is
similar to analogous results in \cite[Lemmas~4 and~5]{DelshamsG04}
and \cite[Lemma~3]{DelshamsG03}.
At first order in $\mu$,
the coefficients of the splitting potential can be approximated,
neglecting the error term in the Melnikov
approximation~(\ref{eq:melniapproxM}),
by the coefficients of the Melnikov potential,
given in~(\ref{eq:alphabeta}):
\ $\Lc_k\sim\mu L_k=\mu\alpha_k\,\ee^{- \beta_k}$.
\ As mentioned in Section~\ref{sect:gn},
the main behavior of the coefficients~$L_k(\eps)$ is given by
the exponents $\beta_k(\eps)$, which have been written in~(\ref{eq:beta_gk})
in terms of the functions $g_k(\eps)$.
In particular, the coefficient~$L_{S_1}$
associated to the dominant harmonic $k=S_1(\eps)$
can be expressed in terms of the function $h_1(\eps)$
introduced in~(\ref{eq:h1}).
In this way, we obtain an estimate for the factor $\ee^{-\beta_{S_i}}$,
which provides the exponential factor in~(a).

We also consider the factor $\alpha_k$,
with $k=S_1(\eps)$. Recalling from~(\ref{eq:estimS}) that
\ $\abs{S_1}\sim\eps^{-1/6}$,
we get from~(\ref{eq:alphabeta}) that
\ $\alpha_{S_1}\sim\eps^{-1/6}$,
\ which provides the polynomial factor in part~(a).

The estimate obtained is valid for the dominant coefficient
of the Melnikov potential $L(\theta)$.
To complete the proof of part~(a),
one has to show that an analogous estimate is also valid
for the splitting potential $\Lc(\theta)$,
i.e.~when the error term in the Poincar\'e--Melnikov
approximation~(\ref{eq:melniapproxM}) is not neglected.
This requires to obtain an upper bound
(provided in \cite[Th.~10]{DelshamsGS04})
for the corresponding coefficient of the error term in~(\ref{eq:melniapproxM})
and show that, in our singular case $\mu=\eps^r$, it is also exponentially
small and dominated by the main term in the approximation.
This can be worked out straightforwardly as in \cite[Lemma~5]{DelshamsG04}
(where the case of the golden number was considered),
so we omit the details here.

The proof of part~(b) is carried out in similar terms.
For the dominant harmonic $k=S_2(\eps)$ inside the set $\Zc_2$,
we also get $\abs{S_2}\sim\eps^{-1/6}$ as in~(\ref{eq:estimS}),
and an exponentially small estimate for $\Lc_{S_2}$
with the function $h_2(\eps)$ defined in~(\ref{eq:h2}).
Such estimates are also valid if one considers the whole sum in~(b),
since for any given $\eps$ the terms of this sum can be bounded
by a geometric series and, hence, it can be estimated by its dominant term
(see \cite[Lemma~4]{DelshamsG04} for more details).
\qed

With regard to the proof of Theorem~\ref{thm:main},
we need to measure the size of the perturbation $\F^{(2)}(\theta)$
in~(\ref{eq:L1F2}) with respect to the coefficient $\Lc_{S_1}$
of the approximation $\Lc^{(1)}(\theta)$.
Since by Lemma~\ref{lm:dominantsL} the size of $\F^{(2)}(\theta)$
is given by the size of its dominant harmonic,
we introduce the following small parameter,
\[
  \eta_{2,1}:=\frac{\Lc_{S_2}}{\Lc_{S_1}}
  \sim\exp\pp{-\frac{C_0(h_2(\eps)-h_1(\eps))}{\eps^{1/6}}},
\]
as a measure of the perturbation $\F^{(2)}(\theta)$ in~(\ref{eq:L1F2}),
relatively to the size of the dominant coefficient $\Lc_{S_1}$.
Although we define the parameter $\eta_{2,1}$
in terms of the coefficients of~$\Lc(\theta)$,
we can also define it from the coefficients of its derivative,
the splitting function~$\M(\theta)=\nabla\Lc(\theta)$,
in view of~(\ref{eq:Mk}) and the fact that the respective factors
have the same magnitude: $\abs{S_1}\sim\abs{S_2}\sim\eps^{-1/6}$.

Notice that the parameter $\eta_{2,1}$ is always exponentially small in $\eps$,
provided we exclude some small neighborhoods
of the ``transition values'' $\ceps$,
where $\Lc_{S_1}$ and $\Lc_{S_2}$ have the same magnitude.

\proofof{Theorem~\ref{thm:main}}
Applying Lemma~\ref{lm:dominantsL}, we see that
the coefficient of the dominant harmonic of the splitting function $\M(\theta)$
is greater than the sum of all other harmonics.
More precisely, we have for $\eps\to0$ the estimate
\beq\label{eq:Mdom}
  \max_{\theta\in\T^3}\abs{\M(\theta)}
  =\abs{\M_{S_1}}(1+\Ord(\eta_{2,1}))
  \sim \abs{\M_{S_1}}
  \sim\abs{S_1}\Lc_{S_1},
\eeq
which implies the result,
using the asymptotic estimate~(\ref{eq:estimS}) for $\abs{S_1}$,
and the asymptotic estimate for $\abs{\M_{S_1}}$, in terms of $h_1(\eps)$,
deduced from Lemma~\ref{lm:dominantsL}(a).

Nevertheless, the previous argument does not apply directly
when $\eps$ is close to a transition value $\ceps$
where $h_1$ and $h_2$ coincide,
i.e.~the first and second dominant harmonics have the same magnitude.
Eventually, more than two harmonics
(but a finite number, according to the arguments given in Lemma~\ref{lm:ZW})
might also have the same magnitude and become dominant.
In such cases, the parameter $\eta_{2,1}$ is not exponentially small,
but we can replace the main term in~(\ref{eq:Mdom}) by a finite number
of terms, plus an exponentially small perturbation, and by the properties of
Fourier expansions the maximum value of~$\abs{\M(\theta)}$ can be
compared to any of its dominant harmonics.
\qed

\paragr{Acknowledgments}
We would like to express our sincere thanks to Carles Sim\'o
for useful discussions and remarks on resonances and Diophantine vectors,
and to Bernat Plans for some useful hints on algebraic number theory.
We also acknowledge the use of EIXAM, the UPC Applied Math cluster system for
research computing (\texttt{https://dynamicalsystems.upc.edu/en/computing}),
and in particular Albert Granados for his support in the use of this cluster.
The author MG also thanks the Dep. de Matem\`atiques i Inform\`atica of the
Univ. de Barcelona for their hospitality and support.

\small
\def\noopsort#1{}

\end{document}